\documentclass[arxiv, 12pt]{elsart}

\usepackage{graphicx}
\usepackage{xcolor}
\usepackage{lineno}

\usepackage{amsmath}
\usepackage{amsfonts}
\usepackage{amssymb}
\usepackage{amscd}
\usepackage{bbm}
\usepackage{bm}

\usepackage{epstopdf,epsfig}
\usepackage{caption,booktabs}
\usepackage{parskip}

\allowdisplaybreaks

\makeatletter

\newcommand{\vv}{E_{\varphi}}

\newtheorem{proposition}{Proposition}

\def\diag{{\rm diag }}

\def\v0{ {{\bf 0}} }

\def\vE{ {{\bf E}} }
\def\vf{ {{\bf f}} }
\def\vF{ {{\bf F}} }

\def\vu{ {{\bf u}} }
\def\vT{ {{\bf T}} }

\def\vv{ {{\bf v}} }

\def\vw{ {{\bf w}} }

\def\vx{ {{\textbf x}} }

\journal{Journal of Scientific Computing}

\begin{document}

\begin{frontmatter}


\title{Error Inhibiting Schemes for Initial Boundary Value Heat Equation}
%
%
%
\author{Adi Ditkowski}
\address{School of Mathematical Sciences, Tel Aviv University, Tel Aviv 69978, Israel}
\ead{adid@tauex.tau.ac.il}

\author{Paz Fink Shustin}
\address{School of Mathematical Sciences, Tel Aviv University, Tel Aviv 69978, Israel}
\ead{pazfink@mail.tau.ac.il}


\begin{abstract}
Finite Difference (FD) schemes are widely used in science and engineering for approximating solutions of partial differential equations (PDEs). Error analysis of FD schemes relies on estimating the truncation error at each time step. This approach usually leads to a global error whose order is of the same order of the truncation error. For classical FD schemes the global error is indeed of the same order as the truncation error. A particular class of FD schemes is the Block Finite Difference (BFD) schemes, in which the grid is divided into blocks. The structure of such schemes is similar to the structure of the Discontinuous Galerkin (DG) method \cite{baumann1999discontinuous,zhang2003analysis,shu2009discontinuous}, and allows inhabitation of the truncation errors. Recently, much effort was devoted to design BFD schemes whose global error converges faster than the truncation error (see \cite{ditkowski2015high} for PDEs and \cite{ditkowski2017error} for ODES). 

In this paper, we elaborate the approach presented in \cite{ditkowski2015high} for the heat equation with periodic boundary conditions, and it is a summary of \cite{FinkPaz2017Eisf}. We generalize this methodology to design BFD schemes for the heat equation with Dirichlet or Neumann boundary conditions, whose global error converges faster than the truncation error. Such schemes are henceforth called Error Inhibiting Schemes. We provide an explicit error analysis, including proofs of stability and convergence of the proposed schemes. We illustrate our approach using several numerical examples, which demonstrate the efficiency of our method in comparison to standard FD schemes.
%
\end{abstract}


\end{frontmatter}

\section{Introduction}\label{intro}
Consider a PDE of the form:
\begin{eqnarray}\label{1.10}
\hspace{3cm}\frac{\partial \, u}{\partial t} &=&  P \left ( \frac{\partial \,
}{\partial x} \right ) u, \quad  x\in\Omega\subset
{\mathbb{R}}^d \;,
t\ge0 \nonumber \\
u(x,0) &=& f(x) \ .
\end{eqnarray}
Here  $P \left ( {\partial \,
}/{\partial x} \right ) $ is a linear  differential operator with appropriate boundary conditions. We assume that \eqref{1.10} is well-posed, i.e. there exists $K(t)<\infty$ such that $|| u(t)|| \le K(t) || f ||$. Typically, $K(t) = K e^{\alpha t}$.

The solution of \eqref{1.10} is often approximated via semi-discretization of the spatial operator $P \left ( {\partial \,
}/{\partial x} \right )$, which will be denoted by $Q$. We assume the following:
\begin{itemize}
\item \textbf{Assumption 1:} The operator $Q$ is induced by the grid points $\{x_j \}$, $j=1, \ldots,N$.
\item \textbf{Assumption 2:} There exists a matrix $H$ satisfying $H=H^*, \  mI\leq H\leq MI, \ 0<m<M$ and a scalar product
\begin{equation} \label{1.12}
\left( \vu, \vv\right)_H = \left( \vu, H \vv \right )_h = h \sum_{j=1}^N \bar{u}_j w_j, \quad \vw =  H \vv
\end{equation}
such that $Q$ is semi-bounded with respect to $\left( \cdot, \cdot \right)_H$, i.e.,
\begin{equation*} \label{1.20}
\left( \vu, Q \vu\right )_H \,  \leq   \, \alpha \left( \textbf{u}, \vu
\right )_H \, = \, \alpha \left\| \vu   \right\|_H^2
\end{equation*}
for some $\alpha>0$ and $h>0$ the maximal distance between two grid points. Note that $\left( \cdot, \cdot \right)_H$ is equivalent to the standard Euclidean inner product.
\item \textbf{Assumption 3:} Let the local truncation error of $Q$, be defined as
\begin{equation*} \label{1.30}
 \left(\vT_e\right)_j= \left( P w(x_j) \right)\, - \,  \left(Q \vw \right)_j ,\quad j = 1,\dots, N
\end{equation*}
where $w(x)$ is a smooth function and $\vw$ is the projection of $w(x)$ onto the grid. We assume that $ \left\| \vT_e \right\|_H \xrightarrow{N
\rightarrow \infty} 0$.
\end{itemize}

Consider the semi--discrete approximation:
\begin{eqnarray}\label{1.40}
\hspace{4cm} \dfrac{\partial \, \vv}{\partial t} &=&  Q \vv, \quad
t\ge0 \nonumber \\
\vv(0) &=& \vf
\end{eqnarray}
\begin{proposition}
Under Assumptions 1-3 the semi-discrete approximation \eqref{1.40} converges to the solution of \eqref{1.10}.
\end{proposition}
\textbf{Proof}: Let $\vu$ be the projection of $u(x,t)$ onto the grid. From Assumption 3 we have
\begin{equation}\label{1.50}
\frac{\partial \, \vu}{\partial t} \,=\, P\left(\frac{\partial}{\partial x}\right) \vu \,=\, Q \vu +  \vT_e
\; .
\end{equation}
Denote the approximation error by $\vE = \vu - \vv$. By subtracting \eqref{1.40} from
\eqref{1.50}, one obtains the error equation
\begin{equation}\label{1.60} 
\frac{\partial \, \vE}{\partial t} \,=\,  Q \vE +  \vT_e\; .
\end{equation}
By using Assumption 2 we have
\begin{eqnarray*}\label{1.70}
\frac{1}{2}\frac{\partial \, }{\partial t}  \left(\vE,  \vE  \right)_H &=& \left(\vE,   \frac{\partial \, \vE}{\partial t} \right)_H \overset{\eqref{1.60}}{=}\left(\vE,  Q \vE  \right)_H  + \left(\vE,  \vT_e  \right)_H \\
& \leq & \alpha \left\| \vE \right\|_H \, + \, \left\| \vE
\right\|_H \left\| \vT_e \right\|_H
\end{eqnarray*}
implying
\begin{equation}\label{1.80} 
\frac{\partial }{\partial t}  \left\| \vE \right\|_H \,\le \, \alpha
\left\| \vE \right\|_H +  \left\| \vT_e \right\|_H \; .
\end{equation}
Therefore, by Assumption 3
\begin{equation}\label{1.90} 
\left\| \vE \right\|_H(t) \,\le \,\left\| \vE (0) \right\|_H {\rm
e}^{\alpha  t }  + \frac{e^{\alpha  t }-1}{\alpha } \max_{0
\le\tau \le t} \left\| \vT_e \right\|_H \; \xrightarrow{N \rightarrow
\infty} \; 0\; 
\end{equation}
\hfill $ \square$

Note that we assume $\left\| \vE (0)\right\|_H$ is either 0, or at
least of the order of machine accuracy.

Equation \eqref{1.90} shows that if the scheme is stable and consistent,
the numerical approximation $\vv$ converges to the projection of the
exact solution onto the grid, $\vu$.  Furthermore, \eqref{1.90} guarantees that
the global error is bounded by the truncation error $\left\| \vT_e
\right\|_H$. This is one part of the landmark Lax-Richtmyer
equivalence theorem for semi-discrete approximations \cite{lax1956survey,quarteroni2010numerical,gustafsson2013time}. Typically the error and the truncation error are of the same order \cite{gustafsson2013time,morton2005numerical,strikwerda2004finite,allen2011numerical,isaacson2012analysis}. However, note that \eqref{1.60} is an equality whereas \eqref{1.90} is an upper bound on the error norm.

Our goal is to design BFD schemes whose global error is of higher order than the truncation error. We are motivated by \cite{zhang2003analysis}, where Nodal-basis DG schemes were presented and analyzed as BFD schemes. In this paper it was also pointed out that the increasing of accuracy may be related to a phenomenon called Supra–Convergence \cite{kreiss1986supra}. Although the truncation error is of first order, the resulting scheme is of second order. Such BFD schemes for the heat equation with periodic boundary conditions were presented in \cite{ditkowski2015high}.

Our main contribution lies in considering initial-boundary value problems (IBVPs), where the main challenge for constructing high-order schemes is in the design of boundary stencils. Near the boundary, the internal stencils, which are usually central, cannot be used. The boundary schemes should approximate the PDE and boundary conditions while maintaining stability. There are several approaches for constructing these boundary stencils, such as penalty methods \cite{Gottlieb2007Spectral,hesthaven2000spectral,utku1982boundary}. An alternative approach is to use extrapolations and the boundary values to generate ghost points outside of the computational domain. An example of this implementation is the inverse Lax–Wendroff method \cite{tan2010inverse}. However, it is well-known that the scheme next to the boundaries can be of a lower accuracy order and preserve accuracy \cite{abarbanel2000error,gustafsson1975convergence,gustafsson1981convergence,svard2006order}. In particular, in \cite{gustafsson1981convergence,svard2006order} it was shown that the boundary conditions can be of two order less for parabolic, incompletely parabolic and second order hyperbolic equations. Here, however, we consider low order truncation errors in most or all of the grid points. 
We exploit the algebraic structure of BFD schemes to separate the subspaces in which the solution and the truncation error lie and construct numerical mechanisms that inhibit the accumulation of the truncation errors. We call them Error Inhibiting schemes (EIS).

This paper is organized as follows: description of the scheme for the heat equation under periodic boundary conditions and the corresponding proofs of convergence are 
presented in Section \ref{sec_heat_per}. In Section \ref{chap2} we describe how to impose boundary conditions for several IBVPs and develop BFD schemes of orders $3$ and $5$. We present numerical
simulations which demonstrate the efficiency of our method and support the developed theory.

\section{Two-Point Block Finite Difference Scheme for Heat Equation with Periodic Boundary Conditions} \label{sec_heat_per}
Consider the following heat equation
\begin{eqnarray}\label{2.05}
\hspace{2cm} \dfrac{\partial u}{\partial t} &=&   \dfrac{\partial^2 u
 }{\partial x^2} +F(x,t), \quad x\in [0, 2\pi)\,,\;
t\ge0 \nonumber \\ 
u(x,0) &=& f(x) \nonumber \\
F(x,t) &=& F(x+2\pi,t), \ f(x)=f(x+2\pi)
\end{eqnarray}
with periodic boundary conditions. This PDE was considered in \cite{ditkowski2015high}.
\subsection{\textbf{Third order Scheme for Heat Equation with Periodic Boundary Conditions}}\label{sec_heat_per_3rd}
\begin{figure}[!htb]
  \begin{center}
    \includegraphics[clip,scale=0.8]{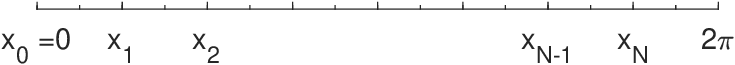}
            \caption{The grid \eqref{2.06}.}\label{fig:grid10}
  \end{center}
\end{figure}

As in \cite{ditkowski2015high}, we introduce a two-point block grid of the following form
\begin{equation}\label{2.06}
x_j = jh, \ x_{j+\frac{1}{2}} = jh+\frac{h}{2}, \quad j=0, \ldots N , \quad h=\frac{2 \pi}{N+1}.
\end{equation}
For simplicity, we assume that $N$ is even. See Figure \ref{fig:grid10}.

Thus, the interval $[0,2\pi]$ is divided into $N+1$ blocks of size $h$, $[x_j,x_{j+1}]$ with a central node $x_{j+1/2}$ where $j = 0,...,N$. Each block contributes two nodes, $x_j,x_{j+1/2}$, to the following  scheme:
\begin{eqnarray}\label{2.10}
\frac{d^2 u_{j}}{dx^2} & \approx & \frac{1}{(h/2)^2} \left[ \left(
u_{j-1/2} -2 u_{j} + u_{j+1/2}  \right) \right . \nonumber \\
			& + &  \left.
					 c \left( -u_{j-1/2} +3 u_{j} - 3 u_{j+1/2} +  u_{j+1} \right)  \right ] \nonumber \\
\frac{d^2 u_{{j+1/2}}}{dx^2}  & \approx & \frac{1}{(h/2)^2} \left[
\left( u_{j} -2 u_{j+1/2} + u_{j+1}  \right) \right . \nonumber \\
			& + & \left. c \left( u_{j-1/2} -3 u_{j} + 3 u_{j+1/2} -  u_{j+1} \right)  \right ]
\end{eqnarray}
with local truncation errors 
\begin{eqnarray*}
\left( T_e \right)_{j} &=& \frac{1}{12} \left( \frac{h}{2}\right)^2 \frac{ \partial^4 \, u_j}{\partial x^4} \\ & + & c \left[ \left( \frac{h}{2}\right) \frac{ \partial^3 \, u_j}{\partial x^3}   + \frac{1}{2} \left( \frac{h}{2}\right)^2 \frac{ \partial^4 \, u_j}{\partial x^4} 
 + \frac{1}{4} \left( \frac{h}{2}\right)^3 \frac{ \partial^5 \, u_j}{\partial x^5}\right]+ O(h^4) \nonumber \\
\left( T_e \right)_{j+\frac{1}{2}} &=& \frac{1}{12} \left( \frac{h}{2}\right)^2 \frac{ \partial^4 \, u_{j+\frac{1}{2}} }{\partial x^4} \nonumber \\
	& + & c \left[ -\left( \frac{h}{2}\right) \frac{ \partial^3 \, u_{j+\frac{1}{2}} }{\partial x^3}   + \frac{1}{2} \left( \frac{h}{2}\right)^2 \frac{ \partial^4 \, u_{j+\frac{1}{2}} }{\partial x^4} 
 - \frac{1}{4} \left( \frac{h}{2}\right)^3 \frac{ \partial^5 \, u_{j+\frac{1}{2}} }{\partial x^5}
\right]  O(h^4).
\end{eqnarray*}
Note that both truncation errors are of order $O(h)$. This scheme was first introduced in \cite{ditkowski2015high} without an explicit error analysis. Here we present the complete analysis, for the first time.

We present \eqref{2.10} in matrix form with non-homogeneous term
\begin{eqnarray}\label{2.25}
\hspace{3.5cm} \dfrac{\partial \, \vv}{\partial t} &=&  Q \vv +\vF(t), \qquad
t\ge0 \nonumber \\
\vv(0) &=& \vf
\end{eqnarray}
where $\vv$ approximates the solution of \eqref{2.05} and $\vF(t)$ and $\vf$ are the projections of $F(x,t)$ and $f(x)$ onto the grid, respectively. We assume that the non-homogeneous term is bounded thus, by Duhamel's principle \cite[Theorem 4.7.2]{gustafsson2013time}, does not affect the stability of the scheme. Therefore, the analysis is done for the homogeneous problem only.

In order to analyze the scheme \eqref{2.10}, it is necessary to diagonalize $Q$. However, $Q$ is not a circulant matrix. Therefore it may not be diagonalized by a DFT matrix. To address this problem, we begin by separating the Fourier spectrum into low and high-frequency modes.  Fix $\omega \in \{-N/2,\ldots,N/2\}$ and let
\begin{equation}\label{2.30}
\nu =\nu(\omega)=\begin{cases}
\omega -(N+1) & \quad \omega >0\\
\omega +(N+1) & \quad \omega \leq0
\end{cases} \ .
\end{equation}
Then, the following identities hold
\begin{equation}\label{2.32} 
 e^{i \omega x_j} = e^{i \nu x_j}, \qquad e^{i \omega
x_{j+1/2}} = -e^{i \nu x_{j+1/2}}.
\end{equation}  
Denote
\begin{equation*} \label{2.34}
e^{i \omega \vx} = \scriptsize\begin{pmatrix}
\vdots\\
e^{i \omega x_j} \\
e^{i \omega x_{j+1/2}} \\
\vdots
\end{pmatrix}, \ e^{i \nu \vx} = \scriptsize \begin{pmatrix}
\vdots \\
e^{i \nu x_j} \\
e^{i \nu x_{j+1/2}} \\
\vdots
\end{pmatrix}
\end{equation*}
and
\begin{equation*} 
M = \operatorname{diag}\left(\mu_1,\mu_2 \right), \ \Sigma = \operatorname{diag}\left(\sigma_1,\sigma_2\right)
\end{equation*} 
with 
\begin{eqnarray*}
\mu_1  &=&  -\frac{4 \sin ^2\left(\frac{\omega \, h}{4}\right)}{\left( h/2 \right)^2}-\frac{8 \, i \, 
   c \, e^{\left(\frac{i\, \omega \, h}{4}\right) } \sin ^3\left(\frac{\omega \, h}{4}\right)}{\left( h/2 \right)^2}, \\
   \mu_2 &=& -\frac{4 \sin ^2\left(\frac{\omega \, h}{4}\right)}{\left( h/2 \right)^2}+\frac{8 \, i \, 
   c \, e^{\left(-\frac{i\, \omega \, h}{4}\right) } \sin ^3\left(\frac{\omega \, h}{4}\right)}{\left( h/2 \right)^2},\\
\sigma_1  &=& -\frac{4 \sin ^2\left(\frac{\nu \, h}{4}\right)}{\left( h/2 \right)^2}-\frac{8 \, i \, 
   c \, e^{\left(\frac{i\, \nu \, h}{4}\right) } \sin ^3\left(\frac{\nu \, h}{4}\right)}{\left( h/2 \right)^2}, \\
    \sigma_2  &=&  -\frac{4 \sin ^2\left(\frac{\nu \, h}{4}\right)}{\left( h/2 \right)^2}+\frac{8 \, i \, 
   c \, e^{\left(-\frac{i\, \nu \, h}{4}\right) } \sin ^3\left(\frac{\nu \, h}{4}\right)}{\left( h/2 \right)^2}.
\end{eqnarray*}
%
By \eqref{2.32}, it can be easily verified that
\begin{equation}\label{2.36}
Q e^{i \omega \vx} = \operatorname{diag}\left(M,\dots,M \right)e^{i \omega \vx}, \quad Q e^{i \nu \vx} = \operatorname{diag}\left(\Sigma,\dots,\Sigma \right)e^{i \nu \vx}
\end{equation}
where $Q$ is defined in \eqref{2.10}.
We look for eigenvectors of the form:
\begin{equation}\label{2.40}
\psi_k(\omega) = \frac{\alpha_k}{\sqrt{2 \pi}} e^{i \omega \vx} \,+ \, \frac{\beta_k}{\sqrt{2 \pi}}e^{i \nu \vx}
\end{equation}
with normalization 
\begin{equation}\label{2.41}
|\alpha_k|^2+|\beta_k|^2=1,\ k=1,2.
\end{equation}
Note that \eqref{2.40} is valid for $\omega \neq 0$. For $\omega=0$, we have $\psi_1(0)=e^{i \textbf{0} \vx}/ \sqrt{2\pi}, \psi_2(0)=e^{i (\textbf{N+1}) \vx}/ \sqrt{2\pi}$.

For $\omega \neq 0$ , we obtain the coefficients $\alpha_k, \beta_k$ and eigenvalues (symbols) $\hat{Q}_k$ by considering nodes $x_j, x_{j+1/2}$ in \eqref{2.36}. This yields the following equations: 
\begin{eqnarray} \label{2.42}
\mu_1 \dfrac{\alpha_k}{\sqrt{2\pi}}e^{i\omega x_j}+\sigma_1 \dfrac{\beta_k}{\sqrt{2\pi}}e^{i\nu x_j} &=&  \hat{Q}_k(\omega) \left (\dfrac{\alpha_k}{\sqrt{2\pi}}e^{i\omega x_j}+ \dfrac{\beta_k}{\sqrt{2\pi}}e^{i\nu x_j}\right) \nonumber \\
\mu_2 \dfrac{\alpha_k}{\sqrt{2\pi}}e^{i\omega x_{j+1/2}}+\sigma_2 \dfrac{\beta_k}{\sqrt{2\pi}}e^{i\nu x_{j+1/2}} &=& \hat{Q}_k(\omega) \left (\dfrac{\alpha_k}{\sqrt{2\pi}}e^{i\omega x_{j+1/2}}+ \dfrac{\beta_k}{\sqrt{2\pi}}e^{i\nu x_{j+1/2}}\right) . 
\nonumber \\
\end{eqnarray}
Denoting $r_k=\beta_k/\alpha_k$ and using the relations \eqref{2.32}  gives
\begin{eqnarray}\label{2.42.10}
\hspace{3cm} \mu_1 + \sigma_1 r_k &=& \hat{Q}_k(\omega) (1+r_k) \nonumber \\
\mu_2 - \sigma_2 r_k  &=& \hat{Q}_k(\omega) (1-r_k) \ .
\end{eqnarray}
Thus, for $\omega \neq 0$ we obtain
\begin{eqnarray} \label{2.46}
\hspace{2cm} r_{1,2}(\omega) &=& \frac{(4-8c)\cos\left( \frac{\omega h}{2} \right) \pm \Delta}{2c \left( 2\sin\left( \frac{\omega h}{2} \right) + \sin \left( \omega h \right) \right)}\,i  \nonumber \\
\hat{Q}_{1,2}(\omega) &=& \frac{-4 +2 c \left( \cos \left( \omega h \right)+3\right) \pm  \Delta   }{2 \left( \frac{h}{2} \right)^2}
\end{eqnarray}
%
%
where
\begin{equation*} \label{2.48}
\Delta   = \sqrt{2 c^2 \cos \left( 2\omega h \right) + 38 c^2+8 (c-1) (3 c-1)
\cos \left( \omega h \right) -32c + 8} \ .
\end{equation*}
It can be verified that for all $c<1/2$, $r_{1,2}$ are imaginary and $\hat{Q}_{1,2}$ are real and negative. Therefore the scheme is von Neumann stable.
Using the normalization $|\alpha_k|^2+|\beta_k|^2=1$ we choose $\alpha_k$, $\beta_k$ to be
\begin{equation}\label{2.50}
\alpha_1 = \frac{1}{\sqrt{1+\left | r_1 \right |^2}}, \hspace{0.8em}
\beta_1 = \frac{r_1}{\sqrt{1+\left | r_1 \right |^2}}, \hspace{0.8em}
\alpha_2 = \frac{ \left | r_2 \right | }{r_2 \sqrt{1+\left | r_2 \right |^2}}, \hspace{0.8em}
\beta_2 = \frac{\left | r_2 \right | }{ \sqrt{1+\left | r_2 \right |^2}}
\end{equation}
%
Finally, for $\omega=0$ we have $\mu_1=\mu_2=0$ and $\sigma_1=\sigma_2= (-4+8c) / \left( h/2 \right)^2$. Similarly to the case $\omega \neq 0$, we obtain the following equations: 
\begin{eqnarray*}\label{2.53}
\mu_1 e^{i \textbf{0} x_j} &=& \hat{Q}_1(0) \, e^{i \textbf{0} x_j} \nonumber \\
\hspace{3.5cm} \mu_2 e^{i \textbf{0} x_{j+1/2}} &=& \hat{Q}_1 (0)\,  e^{i \textbf{0} x_{j+1/2}}
\end{eqnarray*}
which yields $\hat{Q}_1(0)=0$, and
\begin{eqnarray*}\label{2.531}
\hspace{3.5cm} \sigma_1 e^{i(N+1) x_j} &=& \hat{Q}_2(0) \, e^{i(N+1) x_j} \nonumber \\
\sigma_2 e^{i(N+1) x_{j+1/2}} &=& \hat{Q}_2(0) \,  e^{i(N+1) x_{j+1/2}}
\end{eqnarray*}
which yields $\hat{Q}_2(0)=(-4+8c) / \left( h/2 \right)^2$.

\subsubsection{\textbf{Stability of The Scheme}}
Let 
\begin{equation}\label{2.54}
\Psi= \footnotesize
\begin{pmatrix}
\vdots & \vdots & \vdots & \vdots & \vdots & \vdots & \vdots\\
\psi_1\left( -\frac{N}{2}\right) & \psi_2\left(-\frac{N}{2}\right) & \psi_1\left(-\frac{N}{2}+1\right) & \psi_2\left(-\frac{N}{2}+1\right) & \cdots & \psi_1\left(\frac{N}{2}\right) & \psi_2\left(\frac{N}{2}\right)\\
\vdots & \vdots & \vdots & \vdots & \vdots & \vdots & \vdots
\end{pmatrix}
\end{equation}
and note that $\Psi$ is not unitary. By \eqref{2.42}, we can write $Q=\Psi \hat{Q} \Psi^{-1}$. In order to guarantee stability of the scheme, we show that $\Psi$ is invertible and $|| \Psi ||_{h/2}$ and $|| \Psi^{-1} ||_{h/2}$ are uniformly bounded in $N$. Here $||\cdot||_{h/2}$ is the operator norm corresponding to the scalar product \eqref{1.12} with $h$ replaced by $h/2$.
%
%
%
%

\begin{figure}[!htb]
  \begin{center}
    \includegraphics[clip,scale=0.45]{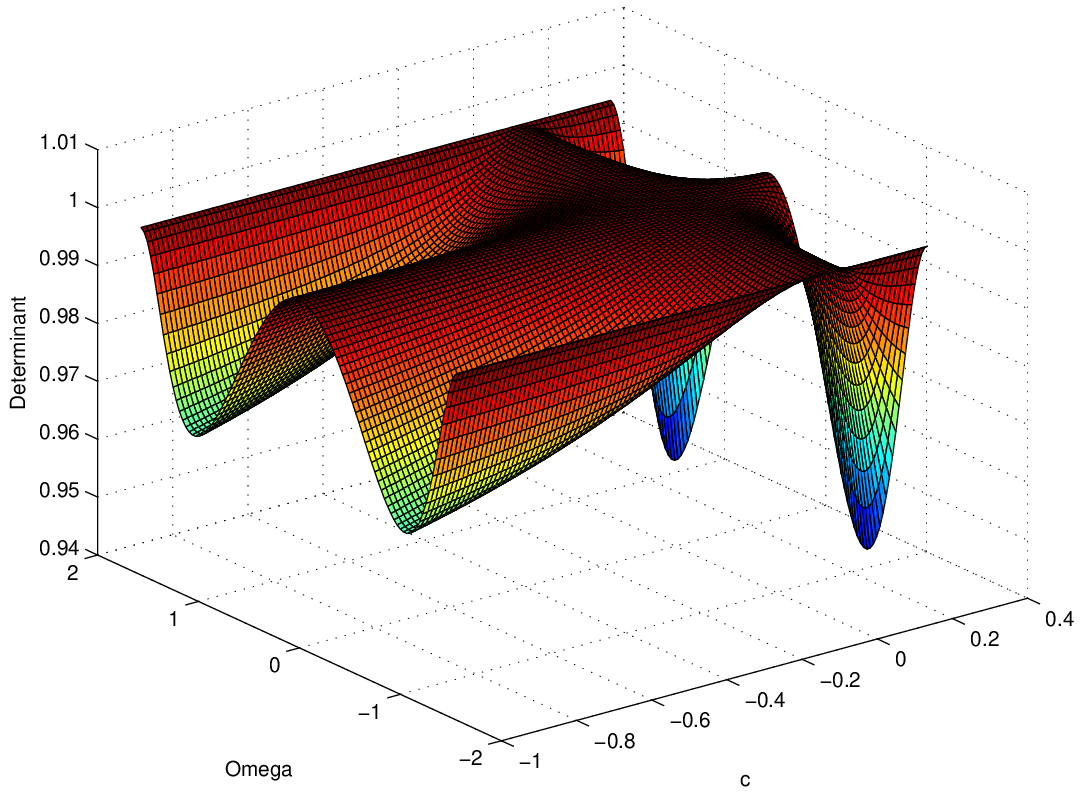}
            \caption{$\operatorname{det}\left(B_{\omega}\right), \quad \omega =-\frac{N}{2},\dots, \frac{N}{2}$}\label{fig:grid11}
  \end{center}
\end{figure}

Denote
\begin{equation*} \label{Psi_norm_10}
{\cal F}^{-1}  = 
\begin{pmatrix}
\vdots & \vdots & \vdots & \vdots & \vdots \\
\frac{e^{i\left (-\frac{N}{2} \right)\vx}}{\sqrt{2\pi} } & \frac{e^{i\left (-\frac{N}{2}+(N+1)\right)\vx}}{\sqrt{2\pi}} & \cdots & \frac{e^{i\left (\frac{N}{2}\right)\vx}}{\sqrt{2\pi}} & \frac{e^{i\left (\frac{N}{2}-(N+1)\right)\vx}}{\sqrt{2\pi}}\\
\vdots & \vdots & \vdots & \vdots & \vdots
\end{pmatrix}.
\end{equation*}
It can be verified that the columns of $\sqrt{\frac{h}{2}} {\cal F}^{-1}$ are orthonormal in the euclidean inner product, which equivalents to 
\begin{equation} \label{ortho_w}
\left( \frac{e^{i\omega_m \textbf{x}}}{\sqrt{2\pi}}, \frac{e^{i\nu_n \textbf{x}}}{\sqrt{2\pi}} \right)_{h/2} = 0, \quad \left( \frac{e^{i\omega_m \textbf{x}}}{\sqrt{2\pi}}, \frac{e^{i\omega_n \textbf{x}}}{\sqrt{2\pi}} \right)_{h/2} = \left( \frac{e^{i\nu_m \textbf{x}}}{\sqrt{2\pi}}, \frac{e^{i\nu_n \textbf{x}}}{\sqrt{2\pi}} \right)_{h/2} = \delta_{n,m} 
\end{equation}
for any $n,m$. For $\omega =-\frac{N}{2},\dots, \frac{N}{2}$ we denote
\begin{equation} \label{B_norm_20}
B_{\omega}= 
\begin{pmatrix}
\alpha_1 \left (\omega \right)& \alpha_2 \left (\omega \right) \\
\beta_1 \left (\omega \right)& \beta_2 \left (\omega \right)
\end{pmatrix}, \quad A = \diag \left( B_{-N/2}, \dots, B_{N/2} \right).
\end{equation}
Then we can write $\Psi = {\cal F}^{-1} A$. Since $\sqrt{\frac{h}{2}} {\cal F}^{-1}$ is unitary, we bound $|| \Psi ||_{h/2}$  by bounding $|| A ||$.
By definition, 
\begin{eqnarray*} \label{B_norm_22}
\Vert B_\omega \Vert &=& \sup_{||y||=1} 
\left\Vert
\begin{pmatrix}
\alpha_1 \left (\omega \right)& \alpha_2 \left (\omega \right) \\
\beta_1 \left (\omega \right)& \beta_2 \left (\omega \right)\\
\end{pmatrix} \cdot 
\begin{pmatrix}
y_1\\
y_2 \\
\end{pmatrix}
\right\Vert \nonumber \\
&=& \sup_{||y||=1} 
\left\Vert
\begin{pmatrix}
\Vert \bm{\alpha} \Vert
\cdot \Vert \bm{y} \Vert \cdot cos(\theta_1)\\
\Vert \bm{\beta} \Vert
\cdot \Vert \bm{y} \Vert \cdot cos(\theta_2)\\
\end{pmatrix}
\right\Vert
\leq \left\Vert
\begin{pmatrix}
\Vert \bm{\alpha} \Vert \\
\Vert \bm{\beta} \Vert \\
\end{pmatrix}
\right\Vert = \sqrt{2}
\end{eqnarray*}
where $\bm{\alpha}=(\alpha_1,\alpha_2),\ \bm{\beta}=(\beta_1,\beta_2),\ \bm{y}=(y_1,y_2) ,\ \theta_1$ is the angle between $\bm{\alpha}$ and $\bm{y}$, and $\theta_2$ is the angle between $\bm{\beta}$ and $\bm{y}$. We also took into account the normalization \eqref{2.41}.
Therefore, given  $\bm{x}\in \mathbb{C}^{2N+2}$ we obtain
\begin{equation*}
\left\Vert  A \bm{x}\right\Vert \leq \left(\left\Vert B_{-N/2} \cdot \left( \begin{matrix}
x_1\\
x_2
\end{matrix} \right) \right\Vert^2 + \ldots + \left\Vert B_{N/2}\cdot \left( \begin{matrix}
x_{2N+1}\\
x_{2N+2}
\end{matrix} \right) \right\Vert^2 \right)^{\frac{1}{2}}=\sqrt{2} \left\Vert \bm{x} \right\Vert
\end{equation*}
which implies $ \left\Vert  A \right\Vert \leq \sqrt{2} $. Therefore, 
\begin{equation*} \label{Psi_norm_40} 
|| \Psi ||_{h/2} = \left\Vert {\cal F}^{-1} \, A \right\Vert_{h/2} = \left\Vert \frac{h}{2}{\cal F}^{-1} \, A \right\Vert = \sqrt{\frac{h}{2}}\left\Vert A \right\Vert \leq \sqrt{h}.
\end{equation*}
We continue by evaluating $\Psi^{-1} = A^{-1}\, {\cal F}$, 
where $A^{-1}$ is a block diagonal matrix and ${\cal F}=\frac{h}{2} \left({\cal F}^{-1}\right)^*$. From \eqref{B_norm_20} we have
\begin{equation} \label{Psi_norm_60}
B_{\omega}^{-1}=\frac{1}{\alpha_1(\omega) \beta_2(\omega)-\alpha_2(\omega) \beta_1(\omega)}
 \left( 
\begin{matrix}
\beta_2(\omega) & -\alpha_2(\omega)\\
-\beta_1(\omega) & \alpha_1(\omega)
\end{matrix}
\right), \quad \omega =-\frac{N}{2},\dots, \frac{N}{2}.
\end{equation}
Moreover, using \eqref{2.50} we get
\begin{equation*} \label{Psi_norm_70} 
\alpha_1(\omega)\beta_2(\omega)-\alpha_2(\omega)\beta_1(\omega) = \frac{\left| r_2(\omega) \right|}{\sqrt{1+\left| r_1(\omega) \right|^2}\sqrt{1+\left| r_2(\omega) \right|^2}}\left( 1-\frac{r_1(\omega)}{r_2(\omega)}\right) 
\end{equation*}
which is estimated to be between $0.9$ and  1 for $c<3/8$ and $-\pi \le \omega\, h\le \pi$ (see Figure \ref{fig:grid11}).
By similar arguments we have
\begin{equation} \label{Psi_norm_80} 
|| \Psi^{-1}  ||_{h/2}=|| A^{-1} {\cal F} ||_{h/2}=\left\Vert A^{-1} \frac{h}{2} \left({\cal F}^{-1}\right)^* \right\Vert_{h/2} = \left\Vert \frac{h}{2} \sqrt{\frac{h}{2}} A^{-1} \right\Vert \le \frac{5h\sqrt{h}}{9\sqrt{2}}.
\end{equation}
In order to prove stability, we expand the numerical solution with respect to $\psi_k(\omega), \ \omega =-\frac{N}{2},\dots, \frac{N}{2}, \ k=1,2 $:
\begin{equation} \label{2.56}
\textbf{v}(t)=\sum_{\omega=-\frac{N}{2}}^{\frac{N}{2}} \hat{v}_1(t,\omega)\psi_1(\omega)+\hat{v}_2(t,\omega)\psi_2(\omega)
\end{equation}
i.e., $\textbf{v}(t) = \Psi\hat{\textbf{v}}(t)$. Substituting \eqref{2.56} into the scheme \eqref{2.25} yields
\begin{equation} \label{2.58}
\sum_{k=1}^2\sum_{\omega=-\frac{N}{2}}^{\frac{N}{2}} \hat{v}_{k}(t,\omega)_t \psi_k(\omega)=\sum_{k=1}^2\sum_{\omega=-\frac{N}{2}}^{\frac{N}{2}} \hat{Q}_k(\omega)\hat{v}_k(t,\omega)\psi_k(\omega).
\end{equation}
By using \eqref{2.54}, equation \eqref{2.58} can be presented in the following matrix form:
\begin{equation}\label{2.60}
\frac{\partial \hat{\vv}}{\partial t} =  \Lambda  \hat{\vv}
\end{equation}
where 
\begin{equation} \label{2.61}
\Lambda = \diag \left( \hat{Q}_1\left(-\frac{N}{2}\right),\, \hat{Q}_2\left(-\frac{N}{2}\right), \, \ldots \, , \, \hat{Q}_1\left(\frac{N}{2}\right), \, \hat{Q}_2\left(\frac{N}{2}\right) \right) \,.
\end{equation}
The system \eqref{2.60} is a first order linear ODE system with the solution $\hat{\textbf{v}}(t) = e^{ \Lambda t} \hat{\vv}(0)$. 
Hence,
\begin{eqnarray*}\label{2.66}
\Vert \textbf{v} \Vert_{h/2}(t) &=& \Vert \Psi\hat{\textbf{v}} \Vert_{h/2}(t) \leq  \Vert \Psi \Vert_{h/2} \cdot \Vert \hat{\textbf{v}} \Vert_{h/2}(t) \leq 
\Vert \Psi \Vert_{h/2} e^{t\cdot\max \limits_{\omega} {\rm Re}(\hat{Q}_j) } \Vert \hat{\textbf{v}}(0) \Vert_{h/2} \\
& \leq &
\Vert \Psi \Vert_{h/2} \Vert \Psi^{-1} \Vert_{h/2} e^{ t\cdot\max \limits_{\omega} {\rm Re}(\hat{Q}_j)} \Vert \textbf{v}(0) \Vert_{h/2}\\
& \leq & \frac{5h^2}{9\sqrt{2}} e^{ t\cdot\max \limits_{\omega} {\rm Re}(\hat{Q}_j)} \Vert \textbf{v}(0) \Vert_{h/2}
\end{eqnarray*}
where in the last inequality we used $\Vert \Psi \Vert_{h/2}  \Vert \Psi^{-1} \Vert_{h/2} \leq  \frac{5h^2}{9\sqrt{2}} $.  Therefore, the scheme is stable. 

\subsubsection{\textbf{Error Estimation}} \label{err_est}
This section demonstrates that our proposed scheme has a global error that is smaller than the truncation error. To do so, we exploit the interaction between the operator $Q$ and the truncation error. The analysis is performed using the eigenvectors expansion.
We assume that the solution is sufficiently smooth, i.e. $u \in C^5\left[ 0,2\pi \right]$. Therefore, the expansion coefficients, denoted by $\hat{u}(t,\omega)$, decay as $\omega^{-6}$.
%

Consider the expansion of \vE in the eigenvectors basis
\begin{equation} \label{3.12}
\textbf{E}(t)=\sum_{\omega=-\frac{N}{2}}^{\frac{N}{2}} \hat{E}_1(t,\omega)\psi_1(\omega)+\hat{E}_2(t,\omega)\psi_2(\omega)
\end{equation}
which equivalent to write $\vE=\Psi \hat{\vE}$. Hence, from the equation of the error \eqref{1.60} we have
\begin{equation} \label{3.13}
\frac{\partial \hat{\vE}}{\partial t} = \Lambda \hat{\vE}+\hat{\textbf{T}}_e \ .
\end{equation}
%
The solution for the error in norm is: 
\begin{eqnarray} \label{3.14}
\Vert \hat{\vE}(t) \Vert_{h/2}  &=& 
	 \left\Vert e^{ \Lambda t} \hat{\vE}(0) + e^{\Lambda t} \int\limits_0^t e^{-\Lambda \tau}{\hat{\textbf{T}}_e d\tau} \right\Vert_{h/2}  \nonumber \\
	 &\leq& \left\Vert e^{ \Lambda t} \hat{\vE}(0) \right\Vert_{h/2} + \left\Vert e^{\Lambda t} \int\limits_0^t e^{-\Lambda \tau}{\hat{\textbf{T}}_e d\tau} \right\Vert_{h/2} \ .
\end{eqnarray}
%
%

Since the initial error is either 0 or at most of the order of machine error, the term $e^{\Lambda t}\hat{\vE}(0)$ can be neglected. In order to continue the estimation of the error in \eqref{3.14}, we need to evaluate $\hat{\textbf{T}}_e$.
Recall that from the scheme \eqref{2.10} we have
\begin{eqnarray}  \label{3.26}
(T_e)_j &=& \frac{1}{12}\left(  \frac{h}{2}\right)  ^2  \frac{ \partial^4 \, u_j }{\partial x^4} + c\left[ \left( \frac{h}{2}\right) (u_j)_{xxx} + \frac{1}{2}\left(   \frac{h}{2}\right)  ^2  \frac{ \partial^4 \, u_j }{\partial x^4} \right] + O(h^3)  \nonumber \\
(T_e)_{j+1/2} &=& \frac{1}{12}\left(  \frac{h}{2}\right)  ^2  \frac{ \partial^4 \, u_{j+\frac{1}{2}} }{\partial x^4} \nonumber \\
 && \hspace{2em} \, + \, c\left[ -\left( \frac{h}{2}\right) (u_{j+1/2})_{xxx} + \frac{1}{2}\left(   \frac{h}{2}\right)  ^2  \frac{ \partial^4 \, u_{j+\frac{1}{2}} }{\partial x^4} \right] + O(h^3) 
\end{eqnarray}
and denote $\vT_e = \vT_h + \textbf{T}_\ell$, where 
\begin{eqnarray}  \label{3.30}
\vT_h  &=&   c\left(  \frac{h}{2}\right) \diag \left( 1,-1,1,-1,...,1,-1 \right) \textbf{u}_{xxx} +O(h^3) 
 \nonumber \\
\vT_\ell  &=&  \frac{6c+1}{12}\left(  \frac{h}{2}\right)^2  \frac{ \partial^4 \, \vu }{\partial x^4} + O(h^4) \ .
\end{eqnarray}
Now, consider the expansion of $\textbf{u}$ in the eigenvectors basis
\begin{eqnarray}  \label{3.32}
\textbf{u}(t) &=&  \sum_{\omega=-\frac{N}{2}}^{\frac{N}{2}} \hat{u}_1(t,\omega)\psi_1(\omega)+\hat{u}_2(t,\omega)\psi_2(\omega) 
 \nonumber \\
  &=& 
 \frac{1}{\sqrt{2\pi}} \sum_{\omega=-\frac{N}{2}}^{\frac{N}{2}} \hat{u}(t,\omega)e^{i\omega \textbf{x}}+\hat{u}(t,\nu)e^{i\nu \textbf{x}}
\end{eqnarray}
where the last equality is the Fourier expansion which holds with
\begin{equation*}
\hat{u}(t,\omega) = \alpha_1 \hat{u}_1(t,\omega) + \alpha_2 \hat{u}_2(t,\omega), \quad \hat{u}(t,\nu) = \beta_1 \hat{u}_1(t,\omega) + \beta_2 \hat{u}_2(t,\omega) \ .
\end{equation*}
Using the smoothness assumption on $\textbf{u}$, replacing $\vu_{xxx},\frac{ \partial^4 \, \vu }{\partial x^4}$ in \eqref{3.30} with the Fourier representation \eqref{3.32} gives
\begin{eqnarray} \label{3.38}
\textbf{T}_h &=& c\left(  \frac{h}{2}\right) \frac{1}{\sqrt{2\pi}}  \sum_{\omega=-\frac{N}{2}}^{\frac{N}{2}} (i\omega)^3\hat{u}(t,\omega) \diag \left( 1,-1,1,-1,...,1,-1 \right) e^{i\omega \textbf{x}} \nonumber \\
&&  \hspace{2em} \, + \, (i\nu)^3\hat{u}(t,\nu) \diag \left( 1,-1,1,-1,...,1,-1 \right) e^{i\nu \textbf{x}} + O(h^3)   \\
\textbf{T}_\ell &=& \frac{6c+1}{12}\left(  \frac{h}{2}\right)^2 \frac{1}{\sqrt{2\pi}}  \sum_{\omega=-\frac{N}{2}}^{\frac{N}{2}} (i\omega)^4\hat{u}(t,\omega)e^{i\omega \textbf{x}}+(i\nu)^4\hat{u}(t,\nu)e^{i\nu \textbf{x}} + O(h^4) \ .  \nonumber
\end{eqnarray}
%
Next, we use the relation between $\hat{\vT}_e$ and ${\vT}_e$: 
\begin{equation*} \label{3.40}
\hat{\vT}_e={\Psi}^{-1} \vT_e = {\Psi}^{-1} \vT_\ell + {\Psi}^{-1} \vT_h = A^{-1} {\cal F} \vT_\ell + A^{-1} {\cal F} \vT_h
\end{equation*}
and treat the last two terms separately. Recall that ${\cal F} = \frac{h}{2} ({\cal F}^{-1})^{*}$.

First, consider the term $A^{-1} {\cal F} \vT_h$ and begin by looking at the m-position in ${\cal F} \vT_h$ which corresponds the frequency $\omega_m$: 
\begin{eqnarray*}  \label{3.42}
\left( {\cal F} \vT_h \right)_m &= & \frac{h}{2} \frac{1}{\sqrt{2\pi}}e^{-i\omega_m \textbf{x}^T}
 \nonumber \\
&& \hspace{1em} \Bigg(  c\left(  \frac{h}{2}\right) \frac{1}{\sqrt{2\pi}}  \sum_{\omega=-\frac{N}{2}}^{\frac{N}{2}} (i\omega)^3\hat{u}(t,\omega) \diag \left( 1,-1,1,-1,...,1,-1 \right) e^{i\omega \textbf{x}} \\
&  &  \hspace{2em}  \, + \, (i\nu)^3\hat{u}(t,\nu) \diag \left( 1,-1,1,-1,...,1,-1 \right) e^{i\nu \textbf{x}} \Bigg) + O(h^3) \\
& = & c \left( \frac{h}{2} \right)^2 \frac{1}{2\pi} e^{-i\omega_m \textbf{x}^T} \sum_{\omega=-\frac{N}{2}}^{\frac{N}{2}} (i\omega)^3\hat{u}(t,\omega)e^{i\nu \textbf{x}} + (i\nu)^3\hat{u}(t,\nu) e^{i\omega \textbf{x}}  + O(h^3) \\
& = & \frac{ch}{2} \left( i\nu(\omega_m) \right) ^3\hat{u} \left( t,\nu(\omega_m) \right) + O(h^3)
\end{eqnarray*}
%
where the first equality is due to \eqref{2.32} and the last equality is due to the orthogonality relations \eqref{ortho_w}.
Similarly, for the corresponding frequency $\nu(\omega_m)$, we obtain 
\begin{equation*}
\left( {\cal F} \vT_h \right)_{m+1} = \frac{ch}{2} (i\omega_m)^3\hat{u}(t,\omega_m) + O(h^3) \ .
\end{equation*}

In order to evaluate $A^{-1} {\cal F} \vT_h$, we use the Taylor expansions of the terms in the matrix ${B_\omega}^{-1}$ \eqref{Psi_norm_60}. To do so, it is sufficient to look at $\omega h \ll 1$:
\begin{eqnarray} \label{coeff_taylor}
&& \alpha_1 = 1-\frac{c^2 }{32 (1-2
   c)^2} \left(\frac{\omega h}{2} \right)^6+O\left( h^7\right), \quad \beta_1 = -\frac{i c  }{4 - 8 c}\left(
\frac{\omega h}{2} \right)^3 + O(h^5) \nonumber \\
&& \alpha_2 =  \frac{i c  }{2 c-1}\left( \frac{\omega h}{2} \right) + O(h^3), \quad \beta_2 = 1 + O(h^2) .
\end{eqnarray} 
%
Since we consider $\omega h \ll 1$, i.e. small values of $\omega$, then $\omega=O(1)$. This implies that $\nu=O(h^{-1})$. From the smoothness assumption, $\hat{u}(\nu)=O(h^6)$, which yields $(i\nu)^3 \hat{u}(\nu) = O(h^3)$.
Hence, we have
\begin{equation*} \label{3.43}
{B_\omega}^{-1}
\begin{pmatrix}
\left( {\cal F} \vT_h \right)_m \\
\left( {\cal F} \vT_h \right)_{m+1}
\end{pmatrix}=
\begin{pmatrix}
O(1)& O(h)\\
O(h^3) & O(1)
\end{pmatrix}
\cdot \frac{ch}{2}  
\begin{pmatrix}
O(h^3)\\
O(1)
\end{pmatrix} = 
\begin{pmatrix}
O(h^2) \\
O(h)
\end{pmatrix} \ .
\end{equation*}

For the term $A^{-1} {\cal F} \vT_\ell$, the m-position which represent the $\omega_m$ frequency in the vector ${\cal F} \vT_\ell$ is 
\begin{eqnarray*} \label{3.44}
\left( {\cal F} \vT_\ell \right)_m &= & \frac{h}{2} \frac{1}{\sqrt{2\pi}}e^{-i\omega_m \textbf{x}^T}
\Bigg( \frac{6c+1}{12}\left(  \frac{h}{2}\right)^2 \frac{1}{\sqrt{2\pi}}  \nonumber \\
&& \hspace{2em} \sum_{\omega=-\frac{N}{2}}^{\frac{N}{2}} (i\omega)^4\hat{u}(t,\omega) e^{i\omega \textbf{x}} + (i\nu)^4\hat{u}(t,\nu)e^{i\nu \textbf{x}} \Bigg) + O(h^4) \\
& = & \frac{6c+1}{12} \left(  \frac{h}{2}\right)^2 (i\omega_m)^4\hat{u}(t,\omega_m) + O(h^4)
\end{eqnarray*}
%
%
%
where the second equality is due to \eqref{ortho_w}.  Similarly, for the corresponding frequency $\nu(\omega_m)$, we obtain 
\begin{equation*}
\left( {\cal F} \vT_\ell \right)_{m+1} = \frac{6c+1}{12} \left(  \frac{h}{2}\right)^2 (i\nu_m)^4\hat{u}(t,\nu_m) + O(h^4) \ .
\end{equation*}
Thus, similar to the case of ${\cal F} \vT_h$, we have
\begin{equation*} \label{3.46}
{B_\omega}^{-1}
\begin{pmatrix}
\left( {\cal F} \vT_\ell \right)_m \\
\left( {\cal F} \vT_\ell \right)_{m+1}
\end{pmatrix}=
\begin{pmatrix}
O(1)& O(h)\\
O(h^3) & O(1)
\end{pmatrix}
\cdot \frac{6c+1}{12} \left(  \frac{h}{2}\right)^2 
\begin{pmatrix}
O(1)\\
O(h^2)
\end{pmatrix} =
\begin{pmatrix}
O(h^2) \\
O(h^4)
\end{pmatrix} \ .
\end{equation*}

Eventually, we have
\begin{equation*}
\hat{\vT}_h = A^{-1} {\cal F} \vT_h = 
\begin{pmatrix}
\vdots \\
\left( \hat{T}_{h} \right)_j \\
\left( \hat{T}_{h} \right)_{j+\frac{1}{2}} \\
\vdots
\end{pmatrix} = 
\begin{pmatrix}
\vdots \\
O(h^2) \\
O(h) \\
\vdots
\end{pmatrix} , 
\end{equation*}
and
\begin{equation*} 
 \hat{\vT}_\ell = A^{-1} {\cal F} \vT_\ell = 
\begin{pmatrix}
\vdots \\
\left( \hat{T}_{\ell}  \right)_j \\
\left( \hat{T}_{\ell} \right)_{j+\frac{1}{2}} \\
\vdots
\end{pmatrix} = 
\begin{pmatrix}
\vdots \\
O(h^2) \\
O(h^4) \\
\vdots
\end{pmatrix} \ .
\end{equation*}
It can be seen that $\hat{\textbf{T}}_h$ contributes the high frequency terms $\left( \hat{T}_e \right)_{j+1/2}$.
%
Now, the second term on the right-hand-side of \eqref{3.14} can be bounded by
\begin{equation} \label{3.24}
\left\Vert e^{\Lambda t}
\int_0^t e^{-\Lambda \tau} \hat{\textbf{T}}_e d\tau \right\Vert_{h/2}
\leq 
\left\Vert \Lambda^{-1} \left( e^{\Lambda t} - {\bf{I}}_{2N+2} \right) \max\limits_{0\leq \tau\leq t}  \left|  \hat{\textbf{T}}_e \left(\tau\right) \right|  \;\right\Vert_{h/2}
\end{equation}
%
where by $\max\limits_{0\leq \tau\leq t}  \left|  \hat{\textbf{T}}_e \left(\tau\right) \right|$ we mean taking maximum for each term separately in $\hat{\textbf{T}}_e$.
%
Recall that $\Lambda$ has the eigenvalues $\hat{Q}_1(\omega)$ and $\hat{Q}_2(\omega)$ on its diagonal \eqref{2.61}. For $\omega h \ll 1$ we have:
\begin{eqnarray*}
\hat{Q}_1(\omega) &=& -\omega ^2 +\frac{(1+4 c)  \omega
   ^4}{12-24 c} \left( \frac{h}{2} \right)^2 + O(h^4) = O(1) \\
\hat{Q}_2(\omega) &=& -\frac{4 -8 c}{\text{(h/2)}^2}+(1-4 c) \omega ^2 + O(h^2) = O(h^{-2}) \ .
\end{eqnarray*} 
Since $\left( \hat{T}_e \right)_{j+1/2} = O(h)$ and $\hat{Q}_2 \left(\omega\right)$ are negative, for all $\omega =-\frac{N}{2},\dots, \frac{N}{2}$ we have
\begin{equation*}
\frac{e^{\hat{Q}_2 \left(\omega\right)t} -1}{\hat{Q}_2(\omega)} \left| \max\limits_{0\leq \tau\leq t} \left( \hat{T}_e \right)_{j+1/2} \left(\tau\right) \right| = O(h^3) \ .
\end{equation*}
%
%
In addition, for any $\omega \neq 0$:
\begin{equation*} \label{3.50}
\left|
\frac{e^{\hat{Q}_1(\omega)t} -1}{\hat{Q}_1(\omega)} \right| = \left| \frac{e^{(-\omega^2+H.O.T.)t} -1}{-\omega^2+H.O.T.} \right| \leq \frac{1}{\omega^2} \leq  O(1) \ .
\end{equation*}
For $\omega = 0$, recall that $\hat{Q}_1(0)=0$. Thus, from L'Hopital's rule
\begin{equation*} \label{3.51}
\left|
\lim_{s \to 0} \frac{e^{st} -1}{s} \right| = \left| \lim_{s \to 0} \frac{te^{st}}{1} \right| = t = O(1) \ .
\end{equation*}
Since $\left( \hat{T}_e \right)_j = O(h^2)$ , we have 
\begin{equation*}
\frac{e^{\hat{Q}_1 \left(\omega\right)t} -1}{\hat{Q}_1(\omega)} \left| \max\limits_{0\leq \tau\leq t} \left( \hat{T}_e \right)_j \left(\tau\right) \right| = O(h^2)
\end{equation*}
Finally, from \eqref{3.24} we obtain
\begin{equation} \label{3.52}
\Vert \vE \Vert_{h/2} = \Vert \Psi \hat{\vE} \Vert_{h/2} \leq \Vert \Psi \Vert_{h/2} \Vert \hat{\textbf{E}} \Vert_{h/2} \leq \sqrt{h}e^{\Lambda t} \Vert \hat{\textbf{E}} \Vert_{h/2} (0) + O(h^2) \ .
\end{equation}

Note that from the structure of the truncation error \eqref{3.26} we would expect that for $c=-\frac{1}{6}$ the order of the global error also will be higher than 2. However this is not the case, since in \cite{ditkowski2015high} it was shown that the explicit form of the error is:
\begin{eqnarray*}  \label{3.54}
\vE(t)  &=&  \sum_{\omega=-\frac{N}{2}}^{\frac{N}{2}} 
{\rm e}^{-\omega^2 t} \Bigg[ \left({
\frac{(1+4 c) \omega^2 t  }{12-24 c}\left( \frac{\omega h}{2}
\right)^2 + O(h^4)}\right ){\rm e}^{i \omega \textbf{x}} \nonumber \\
 && \hspace{4em} \, + \, \left (
\frac{i c }{4-8 c}\left( \frac{\omega h}{2} \right)^3 + O(h^5)
\right ){\rm e}^{i \nu \textbf{x}} \Bigg ] 
\end{eqnarray*}
%
and it can be seen that $c=-\frac{1}{4}$ yields third-order scheme. This is illustrated in the following  numerical example.
\subsubsection{Numerical Example}
This section shows a concrete numerical example that supports our error analysis of the scheme \eqref{2.10}.
We consider the problem \eqref{2.05}, with initial condition $f(x)$ and homogeneous term $F(x,t)$ such that the solution is $u(x,t)=e^{\cos(x-t)}$.
The scheme \eqref{2.10} was run for $N=32,64,128,256,512$ grid points with fourth order explicit Runga-Kutta time propagator. In Figure \ref{fig:grid1} below, we compare the scheme for different values of $c$ and it can be clearly seen that for the value $c={1}/{4}$, the scheme indeed becomes of third order  at time $t=2\pi$. 
%
%
\begin{figure}[!htb]
  \begin{center}
    \includegraphics[clip,scale=0.45]{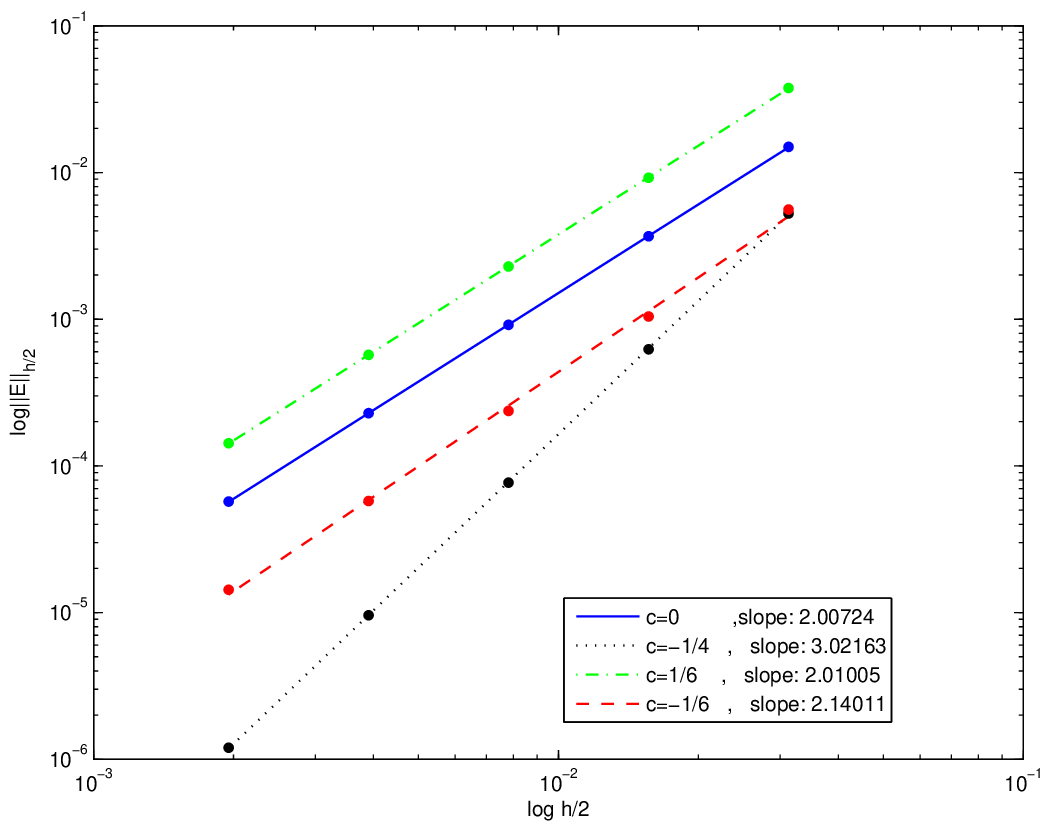}
            \caption{Convergence plot of third Order Scheme for periodic heat equation \eqref{2.10}, $\log_{10} \Vert \textbf{E} \Vert \; vs. \; \log_{10} \left(  \frac{h}{2} \right)  $ for c=0, -1/4, 1/6, -1/6.}\label{fig:grid1}
  \end{center}
\end{figure}

\smallskip
\subsection{\textbf{Fifth Order Scheme: Numerical Example}}
Similarly to the third order scheme \eqref{2.10}, a fifth order scheme can be constructed by using more terms in Taylor expansion, as was presented in \cite{ditkowski2015high}. Consider the following scheme
\begin{eqnarray}  \label{3.56}
\frac{d^2}{dx^2}u_{j} & \approx & \frac{1}{12(h/2)^2} \left[ (-u_{j-1} + 16u_{j-1/2} - 30u_{j} + 16u_{j+1/2} - u_{j+1}) \right.  \nonumber \\
&  &  \left.  \, + \, c(-u_{j-1} + 5u_{j-1/2} - 10u_{j} + 10u_{j+1/2} - 5u_{j+1} + u_{j+3/2}) \right] 
 \nonumber \\
\frac{d^2}{dx^2}u_{j+1/2} & \approx & \frac{1}{12(h/2)^2} \left[(-u_{j-1/2} + 16u_{j} - 30u_{j+1/2} + 16u_{j+1} - u_{j+3/2}) \right.  \nonumber \\
&  & \left. \, + \, c(u_{j-1} - 5u_{j-1/2} + 10u_{j} - 10u_{j+1/2} + 5u_{j+1} - u_{j+3/2}) \right] \ .
 \nonumber \\
\end{eqnarray}
%
%
This is a fourth-order scheme that becomes fifth order for $c={4}/{13}$.
We run this scheme for the same problem from the previous example, for $N=32,64,128,256,512$ grid points and fourth-order explicit Runge-Kutta time propagator. 
In Figure \ref{fig:grid2} below, we compare the scheme for different values of $c$ and it can be clearly seen that for the value $c={4}/{13}$ the scheme indeed becomes of fifth order at time $t=2\pi$. 

\begin{figure}[!htb]
  \begin{center}
    \includegraphics[clip,scale=0.5]{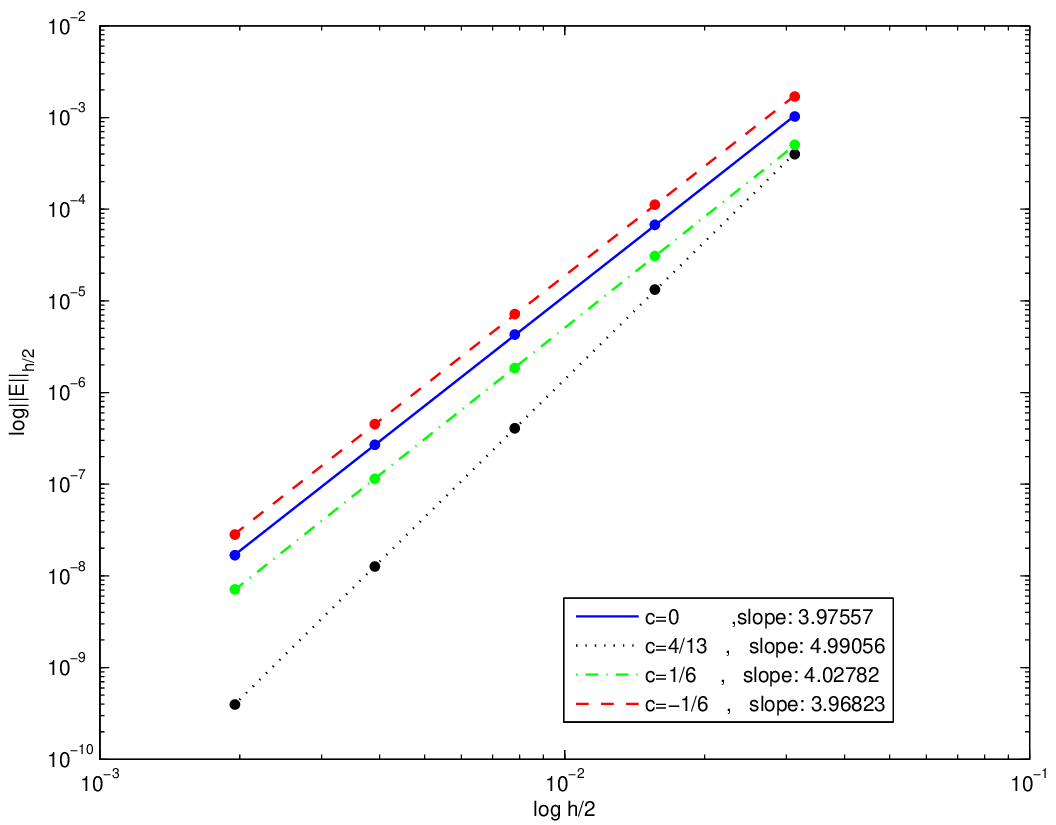}
            \caption{Convergence plot of fifth Order Scheme for periodic heat equation \eqref{2.10}, $\log_{10} \Vert \textbf{E} \Vert \; vs. \; \log_{10} \left(  \frac{h}{2} \right)  $ for c=0, 4/13, 1/6, -1/6.}\label{fig:grid2}
  \end{center}
\end{figure}

\section{Two-point Block Finite Difference Schemes for IBVP Heat Equation} \label{chap2}
In this section, we present two-point BFD schemes and their analysis for the heat equation with non-periodic boundaries. The schemes are developed and analyzed for Dirichlet and Neumann boundary problems.

Consider the following IBVP of the non-homogeneous heat equation: 
\begin{eqnarray} \label{4.10}
\hspace{2cm} \frac{{\partial}u}{{\partial}t} &=& \frac{{\partial}^2}{{\partial}x^2}u+F(x,t),\quad x \in (0,\pi),\; t\geq0 \nonumber \\
u(x,0) &=& f(x)
\end{eqnarray}
with Dirichlet boundary conditions: 
\begin{equation}\label{4.12} 
u(0,t)=g_0(t),\; u(\pi,t)=g_\pi(t)
\end{equation}
or Neumann boundary conditions:
\begin{equation} \label{4.14} 
u_x(0,t)=g_0(t),\; u_x(\pi,t)=g_\pi(t) \ .
\end{equation}
\subsection{\textbf{Reformulation of The Third Scheme for the Periodic Problem}}
\begin{figure}[!htb]
  \begin{center}
    \includegraphics[clip,scale=0.8]{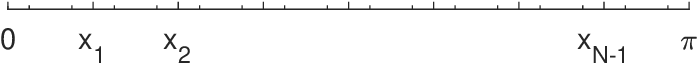}
            \caption{The grid \eqref{4.16} for $[0,\pi]$.}\label{fig:grid3}
  \end{center}
\end{figure}

Although numerical experiments show that the scheme \eqref{2.10} could be generalized for IBVPs using the grid \eqref{2.06}, the analysis becomes much more cumbersome and less intuitive. Therefore, we apply this scheme on a different grid, which will be also applied for the IBVP problem. The new grid does not include the boundaries and allows us to use ghost points for the IBVP problem:
\begin{equation} \label{4.16}
x_{j+1/4} = jh + \frac{h}{4}, \ x_{j+\frac{3}{4}} =  jh + \frac{3h}{4}, \quad j=0,\ldots,N-1, \quad h=\frac{\pi}{N}.
\end{equation}
Note that also $x_N=\pi$. This grid is derived from the grid of the periodic problem on $[0,2\pi]$ where $j=0,\ldots,2N-1$ and $x_{2N}=2\pi$.
For simplicity, we again assume that $N$ is even. See Figure \ref{fig:grid3}.\\

The interval $[0,\pi]$ is divided into $N$ blocks of size $h$, i.e., the interval $[0,2\pi]$ is divided into  $2N$ blocks. Each block contributes two nodes $x_{j+1/4},x_{j+3/4}$ to the following scheme:
\begin{eqnarray}  \label{5.0}
\frac{d^2 }{dx^2} u_{j+1/4} & \approx & \frac{1}{(h/2)^2} \Big[ \left(
u_{j-1/4} -2 u_{j+1/4} + u_{j+3/4}  \right) \nonumber \\
&& \hspace{4em}  \, + \, c \left( -u_{j-1/4} +3 u_{j+1/4} - 3 u_{j+3/4} +  u_{j+5/4} \right)  \Big]
 \nonumber \\
\frac{d^2 }{dx^2} u_{{j+3/4}} &\approx & \frac{1}{(h/2)^2} \Big[
\left( u_{j+1/4} -2 u_{j+3/4} + u_{j+5/4}  \right) \nonumber \\
&&  \hspace{4em}  \, + \, c \left(
u_{j-1/4} -3 u_{j+1/4} + 3 u_{j+3/4} -  u_{j+5/4} \right)  \Big] \ .
\end{eqnarray}
%
For the periodic problem on $[0,2\pi]$, the analysis of this scheme is done in the same manner as the analysis in the previous section.

We begin by splitting the Fourier spectrum into low and high frequency modes. For $\omega=-N+1, \cdots,  N$, let
\begin{equation} \label{5.02}
\nu =\nu(\omega)=\begin{cases}
\omega -2N & \quad \omega >0\\
\omega +2N & \quad \omega \leq0
\end{cases}
\end{equation}
where the following relations  hold:
\begin{eqnarray} \label{ew_4grid}
\forall \omega>0: \ e^{i\nu x_{j+1/4}} &=& -ie^{i\omega x_{j+1/4}}, \quad e^{i\nu x_{j+3/4}}=ie^{i\omega x_{j+3/4}} \nonumber \\ 
\forall \omega\leq 0: \ e^{i\nu x_{j+1/4}} &=& ie^{i\omega x_{j+1/4}}, \quad e^{i\nu x_{j+3/4}}=-ie^{i\omega x_{j+3/4}} \ .
\end{eqnarray}

Note that the orthogonality \eqref{ortho_w} also holds for the grid \eqref{4.16}.
As in the previous section, we look for eigenvectors of the form of \eqref{2.40}. A similar equations to \eqref{2.42.10} can be derived, and the obtained symbols $Q_1(\omega),Q_2(\omega)$ are defined by the same formula in \eqref{2.46}. However, $r_1$ and $r_2$ are different when $\omega \neq 0$:
\begin{equation*} \label{5.01}
\tilde{r}_{1,2}(\omega)=\frac{(-4+8c)\cos(\omega(h/2))\mp \Delta}{2c(2\sin(\omega(h/2))+\sin (2\omega (h/2)))} \ .
\end{equation*}
and satisfy $\tilde{r}_{1,2}(\omega) = i r_{1,2}(\omega)$ for $\omega>0$ and $\tilde{r}_{1,2}(\omega) = -i r_{1,2}(\omega)$ for $\omega<0$. Respectively, the coefficients \eqref{2.50} vary only by a constant. Hence, we consider the eigenvectors 
\begin{equation} \label{5.03} 
\tilde{\psi}_k(\omega) = \frac{\tilde{\alpha}_k}{\sqrt{2\pi}} e^{i \omega \vx} \,+ \,  \frac{\tilde{\beta}_k}{\sqrt{2\pi}} \, e^{i \nu \vx}, \, k=1,2
\end{equation}
instead of the form \eqref{2.40}.
Note that for $\omega=N$ the eigenvectors are 
\begin{equation*}
\frac{1}{\sqrt{2\pi}} \left( e^{i N \vx}+e^{-i N \vx} \right) , \ \frac{1}{\sqrt{2\pi}} \left( e^{i N \vx}-e^{-i N \vx} \right) 
\end{equation*} which can be expressed as $\frac{1}{\sqrt{2\pi}}\cos(N \vx), \frac{1}{\sqrt{2\pi}}\sin(N \vx)$ respectively.\\
For $\omega=0$ the eigenvectors are
\begin{equation*}
\tilde{\psi}_1(0)=\frac{1}{\sqrt{2\pi}}e^{i \textbf{0} \vx}= \frac{1}{\sqrt{2\pi}}\textbf{1}, \tilde{\psi}_2(0)=\frac{1}{\sqrt{2\pi}}e^{2i \textbf{N} \vx}
\end{equation*}
similarly as before.
%
%

\subsection{\textbf{The Eigenvectors of The IBVP}}
The first crucial step for analyzing the IBVP problem is finding the eigenvectors and eigenvalues. To do so, we reflect the interval $[0,\pi]$ to $[0,2\pi]$, and then use the eigenvectors of the periodic  problem on $[0,2\pi]$, as done in the classical analytic problem.
%
In general, denote the eigenvectors as 
\begin{equation} \label{new_eig}
\phi_k(\omega)=\tilde{\psi}_k(\omega)+\tilde{\gamma} \tilde{\psi}_k(-\omega), \ k=1,2 \ .
\end{equation}
Motivated by the classical eigenvectors of the analytic problem, for Dirichlet conditions we consider the eigenvectors $\phi_k(\omega)$ with  $\tilde{\gamma}=-1$ and for Neumann conditions we consider with $\tilde{\gamma}=1$.
%
%
%
Due to the symmetry of the form \eqref{new_eig}, it is sufficient to consider half of the frequencies in \eqref{5.02}, i.e., $\omega=0,\ldots,N$. Indeed, it can be verified that the following hold for every $\omega=1,...,N-1$ and $k=1,2$:
\begin{eqnarray*} \label{5.1}
\hspace{4cm} Q\tilde{\psi}_k(\omega)&=&{\hat{Q}}_k(\omega)\tilde{\psi}_k(\omega) \\
Q\tilde{\psi}_k(-\omega)&=&{\hat{Q}}_k(\omega)\tilde{\psi}_k(-\omega)
\end{eqnarray*}
where $Q$ is the semi-discrete operator of the scheme \eqref{5.0} (same as for the scheme \eqref{2.05}) and the eigenvalues ${\hat{Q}}_k(\omega)$ are the same as in formula \eqref{2.46}, from which it can be seen that $\hat{Q}_k(-\omega)=\hat{Q}_k(\omega)$ . Subtracting or summing these equations yields
\begin{equation*} \label{5.2}
Q\phi_k(\omega)={\hat{Q}}_k(\omega)\phi_k(\omega)
\end{equation*}
with $\tilde{\gamma} = \mp 1$ respectively.
In particular, this means that $\phi_1(\omega),\phi_2(\omega)$ are also eigenvectors of periodic problem.
In addition, note that $\phi_k(\omega)$ with $\tilde{\gamma} = -1$ satisfies the Dirichlet conditions \eqref{4.12}, since for all  $\omega=1,\ldots,N-1$:
\begin{equation*}
\left. \phi_k(\omega) \right|_{x=0,\pi} = \left.  \left[ \tilde{\psi}_k(\omega) - \tilde{\psi}_k(-\omega) \right] \right|_{x=0,\pi} = 0  \ .
\end{equation*}
Hence, they are eigenvectors for the Dirichlet problem. Similarly, $\phi_k(\omega)$ with $\tilde{\gamma} = 1$ satisfies the Neumann conditions \eqref{4.14} and thus they are the eigenvectors for the Neumann problem.

For $\omega=0,N$  with $\tilde{\gamma} = -1$ notice that $\phi_2(0)=\frac{1}{\sqrt{2\pi}}\left( e^{2iN \vx}-e^{-2iN \vx} \right)$ and $\phi_2(N)=\frac{1}{\sqrt{2\pi}} \left( e^{i N \vx}-e^{-i N \vx} \right)$  satisfy Dirichlet conditions \eqref{4.12} and hence they are eqigenvectors for the Dirichlet problem. Similarly, for $\omega=0,N$  with $\tilde{\gamma} = 1$, $\phi_1(0) = \frac{1}{\sqrt{2\pi}}\textbf{1}$ and $\phi_1(N) = \frac{1}{\sqrt{2\pi}} \left( e^{i N \vx}+e^{-i N \vx} \right)$ satisfy Neumann conditions \eqref{4.14} and thus they are eigenvectors for the Neumann problem.
Moreover, in the case of Dirichlet, the eigenvalues can be expressed in terms of normalized sines, and in the case of Neumann, they can be expressed in terms of normalized cosines.

In summary, for the Dirichlet problem the normalized eigenvectors are $\sqrt{\frac{2}{\pi}}\sin(2N \vx)$, $\sqrt{\frac{2}{\pi}} \sin(N \vx)$ and $\sqrt{\frac{2}{\pi}}  \tilde{\alpha}_k \sin(\omega \vx) + \sqrt{\frac{2}{\pi}} \tilde{\beta}_k \sin(\nu \vx)$ where $\omega=1,\ldots,N-1$. For the Neumann problem, the normalized eigenvectors are $\frac{1}{\sqrt{\pi}}\textbf{1},\sqrt{\frac{2}{\pi}} \cos(N \vx)$ and $\sqrt{\frac{2}{\pi}} \tilde{\alpha}_k \cos(\omega x) +\sqrt{\frac{2}{\pi}} \tilde{\beta}_k \cos(\nu \vx)$ for $\omega=1,\ldots,N-1$.
\subsubsection{\textbf{Example of The Eigenvalues for N=6}}
In order to substantiate the results presented above, we present the numeric symbols which obtained in the periodic two-point block approximation with $N=6$. 
%
%
%
In Table \ref{table1} it can be seen that these eigenvalues are splitted between Dirichlet and Neumann problems respectively:

\begin{table}[!htb]
\centering{}
\caption{IBVP Symbols: Dirichlet and Neumann}
\begin{tabular}{c c c c}
\toprule 
$\omega$ & Dirichlet Symbols &  & Neumann Symbols\\
\midrule 
0 & -87.5415 &  & 0\\
1 & -85.5642, \; -0.99994 & & -85.5642, \; -0.99994\\
2 & -79.8974, \; -3.99654 & & -79.8974, \; -3.99654\\
3 & -71.288, \; -8.9584 & & -71.288, \; -8.9584\\
4 & -60.8583, \; -15.7405 & & -60.8583, \; -15.7405\\
5 & -50.1805, \; -23.7481 & & -50.1805, \; -23.7481\\
6 & -29.1805 & & -43.7708\\
\bottomrule 
\end{tabular}
\label{table1}
\end{table} 

\subsection{\textbf{Two-Point Block, Third Order Scheme for Dirichlet IBVP}}
In order to apply the scheme \eqref{5.0} for the grid points near the boundaries $x=0$ and $x=\pi$, we use ghost the  points $x_{-1/4}=-\frac{h}{4}, \ x_{N+1/4}=\pi+\frac{h}{4}$ respectively. At these points, the scheme is computed using extrapolation of two points and the boundaries:
\begin{eqnarray} \label{5.10}
u_{-1/4} &=& -u_{1/4} + 2g_0 + \left( \frac{h}{4}\right)^2 u_{xx}(0,t)+O(h^4) \nonumber \\
u_{N+1/4} &=& -u_{N-1/4} + 2g_\pi+\left( \frac{h}{4}\right) ^2 u_{xx}(\pi,t)+O(h^4)
\end{eqnarray}
where $u_{xx}(0,t)$ and $u_{xx}(\pi,t)$ are expressed using the PDE:
\begin{eqnarray} \label{5.14}
u_{xx}(0,t) &=& u_t(0,t)-F(0,t) = g_0'(t)-F(0,t) \nonumber \\
u_{xx}(\pi,t) &=& u_t(\pi,t)-F(\pi,t) = g_{\pi}'(t)-F(0,t) \ .
\end{eqnarray}
Consider the approximation of two-point block for $j=1,...,N-2$:
\begin{eqnarray} \label{5.18}
\frac{d^2}{dx^2}u_{j+1/4} & \approx & \frac{1}{(h/2)^2} [(u_{j-1/4}-2u_{j+1/4}+u_{j+3/4})
 \nonumber \\
 && \hspace{4em} \, + \,
 c(-u_{j-1/4}+3u_{j+1/4}-3u_{j+3/4}+u_{j+5/4})]  \nonumber \\
\frac{d^2}{dx^2}u_{j+3/4} & \approx & \frac{1}{(h/2)^2} [(u_{j+1/4}-2u_{j+3/4}+u_{j+5/4})
 \nonumber \\
 && \hspace{4em} \, + \, c(u_{j-1/4}-3u_{j+1/4}+3u_{j+3/4}-u_{j+5/4})]
\end{eqnarray}
whereas near the boundaries we have:
\begin{eqnarray} \label{5.20}
\frac{d^2}{dx^2}u_{1/4} & \approx & \frac{1}{(h/2)^2} \Big[ 2(1-c)g_0+(1-c)\left( \frac{h}{4} \right) ^2 u_{xx}(0,t)  \nonumber \\
&& \hspace{4em}  \, + \, (-3u_{1/4}+u_{3/4})+c(4u_{1/4}-3u_{3/4}+u_{5/4}) \Big ]   \nonumber \\
\frac{d^2}{dx^2}u_{3/4} & \approx & \frac{1}{(h/2)^2}\Big[ 2cg_0+c\left(\frac{h}{4}\right) ^2 u_{xx}(0,t)  \nonumber \\
&& \hspace{4em} \, + \, (u_{1/4}-2u_{3/4}+u_{5/4})+c(-4u_{1/4}+3u_{3/4}-u_{5/4})  \Big ]   \nonumber \\
\frac{d^2}{dx^2}u_{N-3/4} & \approx &  \frac{1}{(h/2)^2}\Big[ 2c g_\pi + c \left( \frac{h}{4}\right) ^2 u_{xx}(\pi,t)  \nonumber \\
&& \hspace{4em} \, + \, (u_{N-5/4}-2u_{N-3/4}+u_{N-1/4})  \nonumber \\
&& \hspace{5em}  \, + \, c(-u_{N-5/4}+3u_{N-3/4}-4u_{N-1/4}) \Big ]   \nonumber \\
\frac{d^2}{dx^2}u_{N-1/4} & \approx & \frac{1}{(h/2)^2} \Big[ 2(1-c)g_\pi+(1-c)\left( \frac{h}{4}\right) ^2 u_{xx}(\pi,t)  \nonumber \\
&& \hspace{4em} \, + \, (u_{N-3/4}-3u_{N-1/4})   \nonumber \\
&& \hspace{5em} \, + \, c(u_{N-5/4}-3u_{N-3/4}+4u_{N-1/4}) \Big ] \ .
\end{eqnarray}
%
Respectively, the local truncation errors for $j=1,...,N-2$ are:
\begin{eqnarray} \label{5.22}
\left( T_e \right)_{j+\frac{1}{4} } & =&
\frac{1}{12} \left( \frac{h}{2}\right)^2 \frac{ \partial^4 \, u_{j+\frac{1}{4} }}{\partial x^4} 
+ c \left[ \left( \frac{h}{2}\right) \frac{ \partial^3 \, u_{j+\frac{1}{4} }}{\partial x^3}  \right .
\nonumber \\
&  & \hspace{2em} \, + \, \left .  \frac{1}{2} \left( \frac{h}{2}\right)^2 \frac{ \partial^4 \, u_{j+\frac{1}{4} }}{\partial x^4} 
 + \frac{1}{4} \left( \frac{h}{2}\right)^3 \frac{ \partial^5 \, u_{j+\frac{1}{4} }}{\partial x^5}
\right]+ O(h^4) = O(h)  \nonumber \\
\left( T_e \right)_{j+\frac{3}{4}} & = &
\frac{1}{12} \left( \frac{h}{2}\right)^2 \frac{ \partial^4 \, u_{j+\frac{3}{4}} }{\partial x^4}  
+ c \left[ -\left( \frac{h}{2}\right) \frac{ \partial^3 \, u_{j+\frac{3}{4}} }{\partial x^3}
\right .
 \nonumber \\
& &\hspace{2em} \, + \, \left . \frac{1}{2} \left( \frac{h}{2}\right)^2 \frac{ \partial^4 \, u_{j+\frac{3}{4}} }{\partial x^4} 
 - \frac{1}{4} \left( \frac{h}{2}\right)^3 \frac{ \partial^5 \, u_{j+\frac{3}{4}} }{\partial x^5}
\right] + O(h^4)  = O(h)  \nonumber \\
\end{eqnarray}
%
%
whereas near the boundaries we have:
\begin{eqnarray}  \label{5.24}
(T_e)_{{\frac{1}{4}}} &=& \frac{15}{192}\left( \frac{h}{2}\right)^2 \frac{ \partial^4 \, u_{{\frac{1}{4}}}}{\partial x^4} + \frac{1}{384}\left( \frac{h}{2}\right)^3 \frac{ \partial^5 \, u_{{\frac{1}{4}}}}{\partial x^5} +
c \left[ \frac{h}{2} \frac{ \partial^3 \, u_{{\frac{1}{4}}}}{\partial x^3} \right .
 \nonumber \\
&  & \left . 
 \, + \, \frac{97}{192}\left( \frac{h}{2}\right)^2 \frac{ \partial^4 \, u_{{\frac{1}{4}}}}{\partial x^4} + \frac{95}{384}\left( \frac{h}{2}\right)^3 \frac{ \partial^5 \, u_{{\frac{1}{4}}}}{\partial x^5} \right] + O(h^4)=O(h)  \nonumber \\
(T_e)_{{\frac{3}{4}}} &=& \frac{15}{192}\left( \frac{h}{2}\right)^2 \frac{ \partial^4 \, u_{{\frac{3}{4}}}}{\partial x^4} + \frac{1}{384}\left( \frac{h}{2}\right)^3 \frac{ \partial^5 \, u_{{\frac{3}{4}}}}{\partial x^5}  +
c \left[ -\frac{h}{2} \frac{ \partial^3 \, u_{{\frac{3}{4}}}}{\partial x^3}  \right .
\nonumber \\
&  & \left . \, + \, \frac{97}{192}\left( \frac{h}{2}\right)^2 \frac{ \partial^4 \, u_{{\frac{3}{4}}}}{\partial x^4} - \frac{95}{384}\left( \frac{h}{2}\right)^3 \frac{ \partial^5 \, u_{{\frac{3}{4}}}}{\partial x^5} \right] + O(h^4)=O(h) \nonumber \\
(T_e)_{N-{\frac{3}{4}}} &=& \frac{15}{192}\left( \frac{h}{2}\right)^2 \frac{ \partial^4 \, u_{N-{\frac{3}{4}}}}{\partial x^4} + \frac{1}{384}\left( \frac{h}{2}\right)^3 \frac{ \partial^5 \, u_{N-{\frac{3}{4}}}}{\partial x^5} +
c \left[ -\frac{h}{2} \frac{ \partial^3 \, u_{N-{\frac{3}{4}}}}{\partial x^3} \right .
 \nonumber \\
&  & \left . \, + \, \frac{97}{192}\left( \frac{h}{2}\right)^2 \frac{ \partial^4 \, u_{{\frac{3}{4}}}}{\partial x^4} - \frac{95}{384}\left( \frac{h}{2}\right)^3 \frac{ \partial^5 \, u_{N-{\frac{3}{4}}}}{\partial x^5} \right] + O(h^4)=O(h)  \nonumber \\
(T_e)_{N-{\frac{1}{4}}} &=& \frac{15}{192}\left( \frac{h}{2}\right)^2 \frac{ \partial^4 \, u_{N-{\frac{1}{4}}}}{\partial x^4} + \frac{1}{384}\left( \frac{h}{2}\right)^3 \frac{ \partial^5 \, u_{N-{\frac{1}{4}}}}{\partial x^5} +
c \left[ \frac{h}{2} \frac{ \partial^3 \, u_{N-{\frac{1}{4}}}}{\partial x^3} \right .
 \nonumber \\
&  &  \left . \, + \, \frac{97}{192}\left( \frac{h}{2}\right)^2 \frac{ \partial^4 \, u_{N-{\frac{1}{4}}}}{\partial x^4} + \frac{95}{384}\left( \frac{h}{2}\right)^3 \frac{ \partial^5 \, u_{N-{\frac{1}{4}}}}{\partial x^5} \right] + O(h^4)=O(h) \ .  \nonumber \\
\end{eqnarray}
Notice that the difference between the truncation errors for $j=1,\ldots,N-2$ and the truncation errors near the boundaries is of order of $O(h^2)$.

The scheme \eqref{5.18},\eqref{5.20} can be written in matrix form with non-homogeneous term:
\begin{eqnarray}\label{5.25}
\hspace{3cm} \frac{\partial \, \vv}{\partial t} &= & Q_D \vv + B_D + \vF(t) 
\quad t\ge0 \nonumber \\
\vv(0) &=& \vf
\end{eqnarray}
where $\vv$ approximates the solution of \eqref{4.10}, $Q_D$ is the discretization matrix of scheme \eqref{5.18}, \eqref{5.20} and
\begin{equation} \label{5.2501}
B_D = \frac{1}{(h/2)^2}
\left(
\begin{array}{c}
2(1-c)g_0+(1-c)\left(\frac{h}{4}\right)^2 u_{xx}(0,t) + O(h^4) \\
2cg_0+c\left(\frac{h}{4}\right)^2 u_{xx}(0,t)  + O(h^4) \\
0 \\
\vdots \\
0 \\
2cg_\pi+c\left(\frac{h}{4}\right)^2 u_{xx}(\pi,t)  + O(h^4) \\
2(1-c)g_\pi+(1-c)\left(\frac{h}{4}\right)^2 u_{xx}(\pi,t) + O(h^4)
\end{array}
\right) \ .
\end{equation}
%
Note that the proving stability can be done as in the periodic problem. Since $\phi_k(\omega)$ are linear combinations of $\psi_k(\omega)$ for $k=1,2$, the main difference is manifested in  constants of the used bounds. Thus, the stability of this scheme is preserved and we omit the details of the proof.

\subsubsection{\textbf{Estimation of The Error}}
In this section, we show our approach to generalize the two-point block scheme for IBVPs yields similar results as presented for the periodic case in Section \ref{err_est}. Most of the analysis is done similarly to the analysis in Section \ref{err_est}, where the main difference is due to the changes in the truncation error. 
We assume that $u \in C^5\left[ 0,\pi \right]$. 
%
%

Consider the error $\vE=\Phi_s \hat{\vE}$, where $\Phi_s$ is the eigenvectors matrix of Dirichlet problem. We can write $\Phi_s = {\cal F}_s^{-1} \tilde{A}$ where
\begin{equation}
{\cal F}_s^{-1}= \scriptsize \sqrt{\frac{2}{\pi}}
\begin{pmatrix}
\vdots & \vdots & \vdots & \vdots & \vdots & \vdots & \vdots \\
\sin (2N \vx) & \sin (N \vx) & \sin \left( \vx \right) & \sin \left( (1-2N)\vx \right) & \cdots &  \sin \left( (N-1)\vx \right) & \sin \left( (-1-N)\vx \right) \\
\vdots & \vdots & \vdots & \vdots & \vdots &\vdots & \vdots
\end{pmatrix}
\end{equation}
and, for $\omega = 1,\dots, N-1$:
\begin{equation} \label{B_w_sin} 
\tilde{B}_0 = \small
\begin{pmatrix}
i \beta_1(0)& \beta_2(0) \\
\alpha_1(N)& -i \alpha_2(N)
\end{pmatrix} , \,
\tilde{B}_{\omega}= 
\begin{pmatrix}
\alpha_1 \left (\omega \right)& -i \alpha_2 \left (\omega \right) \\
i \beta_1 \left (\omega \right)& \beta_2 \left (\omega \right)
\end{pmatrix}, \,
\tilde{A} = \diag \left( \tilde{B}_{0}, \dots, \tilde{B}_{N-1} \right).
\end{equation}

As in Section \ref{err_est}, the obtained solution for the error  in norm satisfies: 
\begin{equation} \label{5.30}
\Vert \hat{\textbf{E}} \Vert_{h/2} \leq
\left\Vert e^{\Lambda t} \hat{\vE}_0 \right\Vert_{h/2} + \left\Vert \Lambda^{-1} \left( e^{\Lambda t} - {\bf{I}}_{2N+2} \right) \max\limits_{0\leq \tau\leq t} \left|  \hat{\textbf{T}}_e \left(\tau\right) \right|  \;  \right \Vert_{h/2}
\end{equation}
%
%
where $\Lambda$ is the eigenvalues matrix as before and $\vT_e=\Phi_s \hat{\vT}_e$.
Since the initial error $\vE_0$ is either 0 or, at most of the order of machine error, the term $e^{\Lambda t}\hat{\textbf{E}}_0$ can be neglected as before.

Denote $\textbf{T}_e = \textbf{T}_I + \textbf{T}_B $, where $\textbf{T}_I$ is the truncation error defined by \eqref{5.22} for $j=0,\ldots,N-1$ and $\textbf{T}_B$ is the difference between the truncation error near the boundaries \eqref{5.22} and the truncation error \eqref{5.24} , i.e.
\begin{eqnarray*} \label{5.30.1}
(T_B)_{{\frac{1}{4}}} &=& -\frac{1}{192}\left( \frac{h}{2}\right)^2 \frac{ \partial^4 \, u_{{\frac{1}{4}}}}{\partial x^4} + \frac{1}{384}\left( \frac{h}{2}\right)^3 \frac{ \partial^5 \, u_{{\frac{1}{4}}}}{\partial x^5} \\
&& \, + \, c \left[ \frac{1}{192}\left( \frac{h}{2}\right)^2 \frac{ \partial^4 \, u_{{\frac{1}{4}}}}{\partial x^4} - \frac{1}{384}\left( \frac{h}{2}\right)^3 \frac{ \partial^5 \, u_{{\frac{1}{4}}}}{\partial x^5} \right] + O(h^4)=O(h^2) \\
(T_B)_{{\frac{3}{4}}} &=& -\frac{1}{192}\left( \frac{h}{2}\right)^2 \frac{ \partial^4 \, u_{{\frac{3}{4}}}}{\partial x^4} + \frac{1}{384}\left( \frac{h}{2}\right)^3 \frac{ \partial^5 \, u_{{\frac{3}{4}}}}{\partial x^5} \\
&& \, + \,  c \left[ \frac{1}{192}\left( \frac{h}{2}\right)^2 \frac{ \partial^4 \, u_{{\frac{3}{4}}}}{\partial x^4} + \frac{1}{384}\left( \frac{h}{2}\right)^3 \frac{ \partial^5 \, u_{{\frac{3}{4}}}}{\partial x^5} \right] + O(h^4)=O(h^2) \\
(T_B)_{N-{\frac{3}{4}}} &=& -\frac{1}{192}\left( \frac{h}{2}\right)^2 \frac{ \partial^4 \, u_{N-{\frac{3}{4}}}}{\partial x^4} + \frac{1}{384}\left( \frac{h}{2}\right)^3 \frac{ \partial^5 \, u_{N-{\frac{3}{4}}}}{\partial x^5} \\
&& \, + \, c \left[ \frac{1}{192}\left( \frac{h}{2}\right)^2 \frac{ \partial^4 \, u_{{\frac{3}{4}}}}{\partial x^4} + \frac{1}{384}\left( \frac{h}{2}\right)^3 \frac{ \partial^5 \, u_{N-{\frac{3}{4}}}}{\partial x^5} \right] + O(h^4)=O(h^2) \\ (T_B)_{N-{\frac{1}{4}}} &=& -\frac{1}{192}\left( \frac{h}{2}\right)^2 \frac{ \partial^4 \, u_{N-{\frac{1}{4}}}}{\partial x^4} + \frac{1}{384}\left( \frac{h}{2}\right)^3 \frac{ \partial^5 \, u_{N-{\frac{1}{4}}}}{\partial x^5} \\
&& \, + \, c \left[ \frac{1}{192}\left( \frac{h}{2}\right)^2 \frac{ \partial^4 \, u_{N-{\frac{1}{4}}}}{\partial x^4} - \frac{1}{384}\left( \frac{h}{2}\right)^3 \frac{ \partial^5 \, u_{N-{\frac{1}{4}}}}{\partial x^5} \right] + O(h^4)=O(h^2) 
\end{eqnarray*}
and $\left( T_B \right)_{j+1/4} = \left( T_B \right)_{j+3/4} $ for $j=1,\ldots,N-2$.
Also, denote $\vT_I = \vT_{I_h} + \vT_{I_\ell}$ where 
\begin{eqnarray} \label{5.30.2}
\vT_{I_h} &=& \frac{c}{12}\left(  \frac{h}{2}\right) \diag \left( 1,-1,1,-1,...,1,-1 \right) \textbf{u}_{xxx}  + O(h^3) \nonumber \\
\vT_{I_\ell} &=& \frac{6c+1}{12}\left(  \frac{h}{2}\right)^2  \frac{\partial^4 \vu}{\partial x^4} + O(h^4) \\
\vT_B &=& \frac{c-1}{192}\left(  \frac{h}{2}\right)^2 \frac{\partial^4 \vu}{\partial x^4} + O(h^4) \ .
 \nonumber 
\end{eqnarray}
%
%
The expansion of $\vu$ in the eigenvectors basis on $[0,\pi]$ is equivalent to the Fourier expansion of its anti-symmetric continuation to $[0,2\pi]$, that is a sine series.
For $\vu_{xxx}$, this means we would have a symmetric continuation to $[0,2\pi]$ with Fourier expansion that is a cosine series.
Thus, we expand $\vu_{xxx}$ in the cosine basis on $[0,\pi]$, and similarly $\frac{\partial^4 \vu}{\partial x^4}$ is expanded in the sine eigenvectors basis:
\begin{eqnarray} \label{5.31}
\vu_{xxx}(t) &=& \sqrt{\frac{2}{\pi}}  \sum_{\omega} \hat{u}_3(t,\omega) \cos(\omega \vx)+\hat{u}_3(t,\nu) \cos(\nu \vx) \nonumber  \\
\frac{\partial^4 \vu}{\partial x^4} (t) &=& \sqrt{\frac{2}{\pi}} \sum_{\omega} \hat{u}_4(t,\omega) \sin(\omega \vx)+\hat{u}_4(t,\nu) \sin(\nu \vx) \ .
\end{eqnarray} 
Note that, if the boundaries are homogeneous for every $u^{(k)}, \ k=1,\ldots,4$, then we could differentiate the expansion of $\vu$ in the eigenvectors basis and obtain that $\hat{u}_3(t,\omega) = -\omega^3 \hat{u}(t,\omega), \ \hat{u}_3(t,\nu) = -\nu^3 \hat{u}(t,\nu), \ \hat{u}_4(t,\omega) = \omega^4 \hat{u}(t,\omega)$ and $\hat{u}_4(t,\nu) = \nu^4 \hat{u}(t,\nu)$ as before.

Substituting \eqref{5.31} into \eqref{5.30.2} gives
\begin{eqnarray} \label{5.32}
\vT_{I_h}  &=& \frac{ch}{2} \sqrt{\frac{2}{\pi}}  \sum_{\omega} \hat{u}_3(t,\omega) \diag \left( 1,-1,1,-1,...,1,-1 \right) \cos(\omega \vx) \nonumber \\ 
&+& \hat{u}_3(t,\nu) \diag \left( 1,-1,1,-1,...,1,-1 \right) \cos(\nu \vx) + O(h^3) \nonumber \\
\vT_{I_\ell} &=& \frac{6c+1}{12}\left(  \frac{h}{2}\right)^2 \sqrt{\frac{2}{\pi}} \sum_{\omega} \hat{u}_4(t,\omega) \sin(\omega \vx) + \hat{u}_4(t,\nu) \sin(\nu \vx) + O(h^4) \nonumber \\
\textbf{T}_B &=& \frac{c-1}{192}\left(  \frac{h}{2}\right)^2 \sqrt{\frac{2}{\pi}} \sum_{\omega} \hat{u}_4(t,\omega) \sin(\omega \vx) + \hat{u}_4(t,\nu) \sin(\nu \vx) + O(h^4)
\, . \quad \quad
\end{eqnarray}

Next, using the relation between $\hat{\vT}_e$ and ${\vT}_e$, we write
\begin{equation*} \label{5.34}
\hat{\vT}_e = {\Phi}_s^{-1}  \vT_e = {\Phi}_s^{-1} \vT_{I_h} + {\Phi}_s^{-1} \vT_{I_\ell} + {\Phi}_s^{-1} \vT_B
\end{equation*}
%
and analyze each term separately.
Note that using \eqref{ortho_w}, it can be shown that for any $n,m$:
\begin{eqnarray} \label{ortho_sin}
\left( \sqrt{\frac{2}{\pi}} \sin(\omega_m \vx), \sqrt{\frac{2}{\pi}} \cos(\omega_m \vx) \right)_{h/2} &&  \hspace{-0.6em}  
= \left( \sqrt{\frac{2}{\pi}} \sin(\omega_m \vx), \sqrt{\frac{2}{\pi}} \sin(\nu_n \vx) \right)_{h/2} 
 \hspace{-0.6em}  = 0 \nonumber \\ 
\left( \sqrt{\frac{2}{\pi}} \sin(\omega_m \vx), \sqrt{\frac{2}{\pi}} \sin(\omega_n \vx) \right)_{h/2} &&  \hspace{-0.6em}  
= \left( \sqrt{\frac{2}{\pi}} \sin(\nu_m \vx), \sqrt{\frac{2}{\pi}} \sin(\nu_n \vx) \right)_{h/2} 
 \hspace{-0.6em}  = \delta_{n,m} \nonumber \\
\left( \sqrt{\frac{2}{\pi}} \cos(\omega_m \vx), \sqrt{\frac{2}{\pi}} \cos(\omega_n \vx) \right)_{h/2} && \hspace{-0.6em} 
 =  \left( \sqrt{\frac{2}{\pi}} \cos(\nu_m \vx), \sqrt{\frac{2}{\pi}} \cos(\nu_n \vx) \right)_{h/2} 
  \hspace{-0.6em}  = \delta_{n,m} \nonumber \\
\end{eqnarray}
and from \eqref{ew_4grid} the following relations obtained for any $j=0,\dots,N-1$:
\begin{eqnarray} \label{sin_cos_sign}
\diag \left( 1,-1,1,-1,...,1,-1 \right) \cos(\omega \vx) &=& - \sin(\nu \vx) \nonumber \\
\diag \left( 1,-1,1,-1,...,1,-1 \right) \cos(\nu \vx) &=& \sin(\omega \vx) \ .
\end{eqnarray}
%
%
Now, for the term ${\Phi}_s^{-1} \vT_{I_h}$ and frequency $\omega_m$, we have: 
\begin{eqnarray*}  \label{5.34.1}
&& \hspace{-3em}   \left( \sqrt{\frac{2}{\pi}} \sin(\omega_m \vx), \vT_{I_h} \right)_{h/2}  \nonumber \\
&=& \frac{h}{2} \sqrt{\frac{2}{\pi}} \sin(\omega_m \vx^T) \Bigg( \frac{ch}{2} \sqrt{\frac{2}{\pi}} \sum_{\omega} \hat{u}_3(t,\omega) \diag \left( 1,-1,1,-1,...,1,-1 \right) \cos(\omega \vx) \\ 
&& \hspace{3em}  \, + \, \hat{u}_3(t,\nu) \diag \left( 1,-1,1,-1,...,1,-1 \right) \cos(\nu \vx) \Bigg) + O(h^3) \\
&=&  c \left( \frac{h}{2} \right)^2 \frac{2}{\pi} \sin(\omega_m \vx^T) \left( \sum_{\omega} - \hat{u}_3(t,\omega) \sin(\nu \vx) + \hat{u}_3(t,\nu) \sin(\omega \vx) \right) + O(h^3) \\
&=& \frac{ch}{2} \hat{u}_3 \left( t,\nu(\omega_m) \right) + O(h^4)
\end{eqnarray*}
%
%
Similarly, for the corresponding frequency $\nu(\omega_m)$, we obtain 
\begin{equation*}
\left( \sqrt{\frac{2}{\pi}} \sin(\nu (\omega_m) \vx), \vT_{I_h} \right)_{h/2} = \frac{ch}{2} \hat{u}_3 \left( t,\omega_m \right) + O(h^4) \ .
\end{equation*}
In order to evaluate ${\Phi}_s^{-1} \vT_{I_h}$, we use the Taylor expansions of the terms in the matrix ${\tilde{B}_\omega}^{-1}$ \eqref{B_w_sin}, which are the same expansions as in \eqref{coeff_taylor} up to constants. To do so, it is sufficient to look at $\omega h \ll 1$, such that as in the periodic problem we have $\omega=O(1)$ and $\nu=O(h^{-1})$. Moreover, the assumption $u \in C^5\left[ 0,\pi \right]$ implies that $u_{xxx} \in C^2\left[ 0,\pi \right]$ and $\frac{\partial^4 u}{\partial x^4} \in C^1\left[ 0,\pi \right]$. However, since the problem is not periodic, we have $\hat{u}_3(t,\nu) = O(h)$ and $\hat{u}_4(t,\nu) = O(h)$. Therefore,
\begin{equation*} 
{\tilde{B}_\omega}^{-1}
\begin{pmatrix}
\left( {\cal F}_s \vT_{I_h} \right)_m \\
\left( {\cal F}_s \vT_{I_h} \right)_{m+1}
\end{pmatrix}=
\begin{pmatrix}
O(1)& O(h)\\
O(h^3) & O(1)
\end{pmatrix}
\cdot \frac{ch}{2}  
\begin{pmatrix}
O(h)\\
O(1)
\end{pmatrix} = 
\begin{pmatrix}
O(h^2) \\
O(h)
\end{pmatrix} \ .
\end{equation*}
For the term ${\Phi}_s^{-1} \vT_{I_\ell}$ and frequency $\omega_m$, we also have: 
\begin{eqnarray*}
&& \hspace{-4em} \left( \sqrt{\frac{2}{\pi}} \sin(\omega_m \vx), \vT_{I_\ell} \right)_{h/2} \\
&=& \frac{h}{2} \sqrt{\frac{2}{\pi}} \sin(\omega_m \vx^T) \Bigg( \frac{6c+1}{12}\left(  \frac{h}{2}\right)^2 \sqrt{\frac{2}{\pi}} \sum_{\omega} \hat{u}_4(t,\omega) \sin(\omega \vx) \\
&&\hspace{4em} \, + \,  \hat{u}_4(t,\nu) \sin(\nu \vx) \Bigg ) + O(h^4) \\
&=& \frac{6c+1}{12} \left(  \frac{h}{2}\right)^2 \hat{u}_3 \left( t,\omega_m \right) + O(h^4)
\end{eqnarray*}
%
%
%
Similarly, for the corresponding frequency $\nu(\omega_m)$, we obtain 
\begin{equation*}
\left( \sqrt{\frac{2}{\pi}} \sin(\nu (\omega_m) \vx), \vT_{I_h} \right)_{h/2} = \frac{6c+1}{12} \left(  \frac{h}{2}\right)^2 \hat{u}_4 \left( t,\nu (\omega_m) \right) + O(h^4) \ .
\end{equation*}
Thus, 
\begin{equation*} 
{\tilde{B}_\omega}^{-1}
\begin{pmatrix}
\left( {\cal F}_s \vT_{I_\ell} \right)_m \\
\left( {\cal F}_s \vT_{I_\ell} \right)_{m+1}
\end{pmatrix}=
\begin{pmatrix}
O(1)& O(h)\\
O(h^3) & O(1)
\end{pmatrix}
\cdot \frac{6c+1}{12} \left(  \frac{h}{2}\right)^2  
\begin{pmatrix}
O(1)\\
O(h)
\end{pmatrix} = 
\begin{pmatrix}
O(h^2) \\
O(h^3)
\end{pmatrix} \ .
\end{equation*}
Finally, since the asymptotic behaviour of the term ${\Phi}_s^{-1} \vT_B$ is similar to ${\Phi}_s^{-1} \vT_{I_\ell}$, we have
\begin{equation*} 
{\Phi}_s^{-1} \vT_{I_h} = 
\begin{pmatrix}
\vdots \\
O(h^2) \\
O(h) \\
\vdots
\end{pmatrix}
, \quad
{\Phi}_s^{-1} \vT_{I_\ell}  = {\Phi}_s^{-1} \vT_B = 
\begin{pmatrix}
\vdots \\
O(h^2) \\
O(h^3) \\
\vdots
\end{pmatrix} \ .
\end{equation*}
From this point, the proof is identical to the corresponding part of the proof in the periodic problem. Therefore, we obtain
\begin{equation*}
||\textbf{E}||_{h/2} \leq ||\Phi_s||_{h/2} \cdot ||\hat{\textbf{E}}||_{h/2} \leq ||\Phi_s||_{h/2} \cdot e^{\Lambda t} ||\hat{\textbf{E}}||_{h/2}(0) + O(h^2)
\end{equation*}
where $||\Phi_s||_{h/2}$ is bounded since $||\Psi||_{h/2}$ is bounded.


\subsection{Numerical Example}
%
%
%
In this section, we demonstrate our theory and analysis for  IBVPs by a numerical example.
Consider the problem \eqref{4.10} with Dirichlet boundary conditions \eqref{4.12}. Also consider the same initial function $f(x)$ and non-homogeneous term $F(x,t)$ such that the solution is $u(x,t)=e^{\cos(x-t)}$, as in the previous example. Note that the boundaries functions are non-homogeneous, since $g_0(t)=e^{\cos(t)}, \ g_{\pi}(t)=e^{\cos(\pi-t)}$. The scheme \eqref{5.18},\eqref{5.20} was run for $N=32,64,128,256,512$ grid points with fourth order explicit Runga-Kutta time propagator. In Figure \ref{fig:grid4} below, we compare the scheme for different values of $c$ and it can be clearly seen that for the value $c=-\frac{1}{4}$, the scheme indeed becomes of third order  at time $t=\pi$, as was demonstrated in the periodic case.
\begin{figure}[!htb]
  \begin{center}
    \includegraphics[clip,scale=0.5]{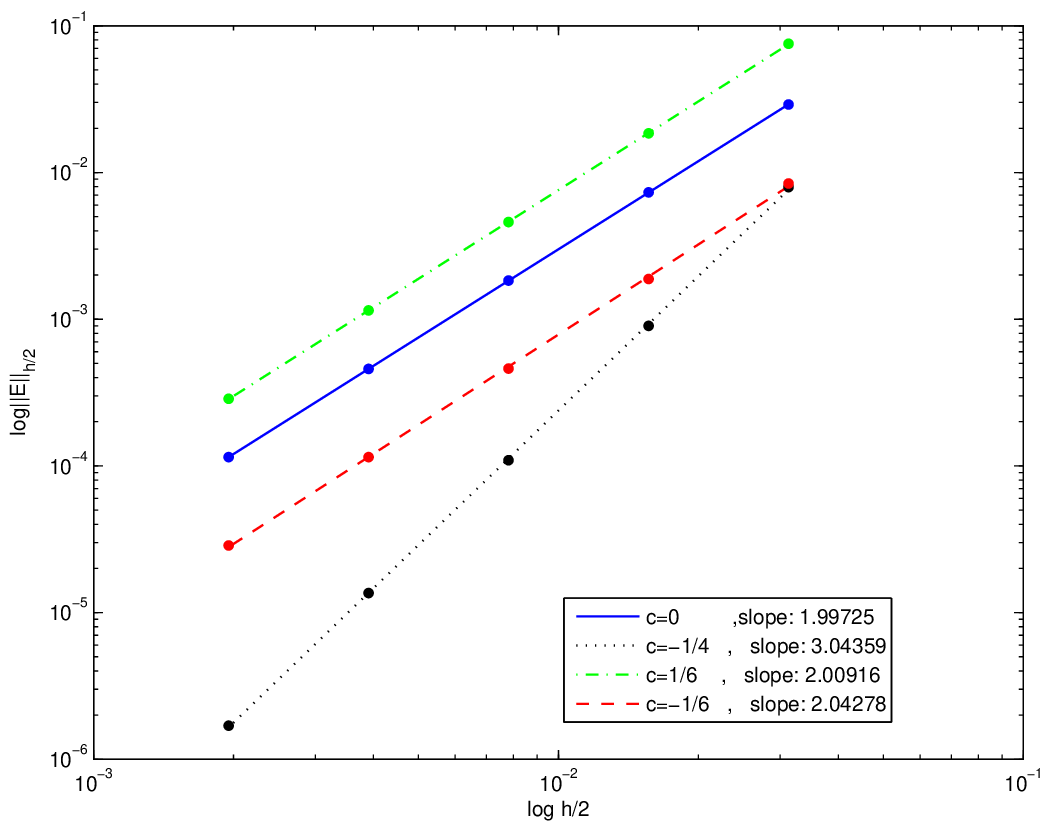}
            \caption{Convergence plot of third order scheme for non-homogeneous Dirichlet \eqref{5.10}, $\log_{10} \Vert \textbf{E} \Vert \; vs. \; \log_{10} \left(  \frac{h}{2} \right)  $ for c=0, -1/4, 1/6, -1/6.}\label{fig:grid4}
  \end{center}
\end{figure}

\subsection{\textbf{Two-Point Block, Third Order Scheme for Neumann IBVP: Numerical Example}}
Consider the problem \eqref{4.10} with Neumann conditions \eqref{4.14}. Since a similar analysis can be derived as was done for the IBVP with Dirichlet conditions, we present only the differences in the boundaries for the scheme \eqref{5.18},\eqref{5.20}. For Neumann conditions, the ghost points are computed using extrapolation of two points and the boundaries:
\begin{eqnarray} \label{5.38}
u_{-1/4} &=& u_{1/4} - \frac{h}{2} u_{x}(0,t) - \frac{1}{3}\left( \frac{h}{4}\right) ^3 u_{xxx}(0,t) + O(h^5) \nonumber \\
u_{N+1/4} &=& u_{N-1/4} + \frac{h}{2} u_{x}(\pi,t) + \frac{1}{3}\left( \frac{h}{4}\right) ^3 u_{xxx}(\pi,t)+O(h^5)
\end{eqnarray}
where $u_{xxx}(0,t)$ and $u_{xxx}(\pi,t)$ are expressed using the PDE:
\begin{eqnarray} \label{5.42}
u_{xxx}(0,t) &=& g_{0}'(t)-F_x(0,t) \nonumber \\
u_{xxx}(\pi,t) &=& g_{\pi}'(t)-F_x(\pi,t)
\end{eqnarray}
%
%
Note that, the vector $B_D$ defined by \eqref{5.2501} is changed respectively and this is the main difference from the case of Dirichlet in the error analysis. 

We demonstrate the case of Neumann conditions using a numerical example. Consider the problem \eqref{4.10} with Neumann boundary conditions \eqref{4.14}. Also consider the same initial function $f(x)$ and non-homogeneous term $F(x,t)$ such that the solution is $u(x,t)=e^{\cos(x-t)}$, as in the previous example. Note that the boundaries functions are $g_0(t)=\sin(t) e^{\cos(t)}, \ g_{\pi}(t)=-\sin(\pi-t) e^{\cos(\pi-t)}$. The scheme \eqref{5.18},\eqref{5.20} was run for $N=32,64,128,256,512$ grid points with fourth order explicit Runga-Kutta time propagator. In Figure \ref{fig:grid5} below, we compare the scheme for different values of $c$ and it can be clearly seen that for the value $c=-\frac{1}{4}$, the scheme indeed becomes of third order  at time $t=\pi$, as was demonstrated in the previous boundary cases.
\begin{figure}[!htb]
  \begin{center}
    \includegraphics[clip,scale=0.5]{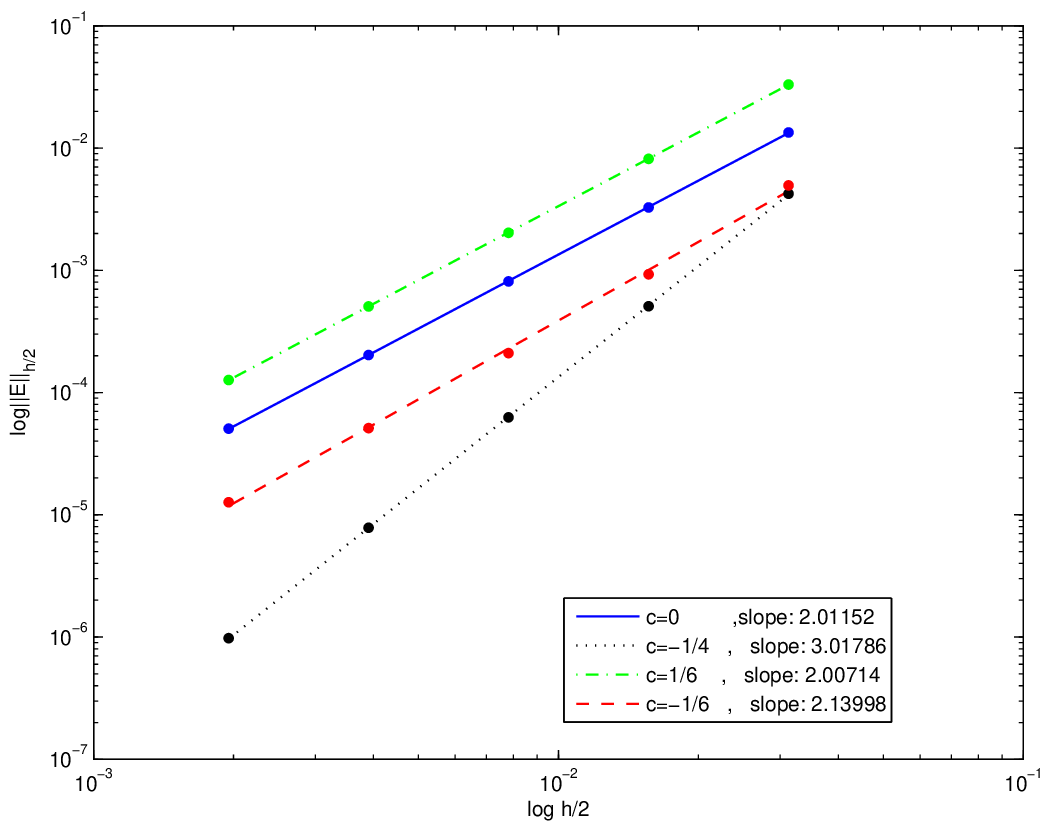}
            \caption{Convergence plot of third order scheme for non-homogeneous Neumann \eqref{4.10},\eqref{4.14}, $\log_{10} \Vert \textbf{E} \Vert \; vs. \; \log_{10} \left(  \frac{h}{2} \right)  $ for c=0, -1/4, 1/6, -1/6.}\label{fig:grid5}
  \end{center}
\end{figure}

\subsection{\textbf{Two-Point Block, Fifth Order Scheme for Dirichlet IBVP: Numerical Example}}
Using our approach, the fifth order scheme \eqref{3.56} that was shown the periodic case can be also adapted to IBVPs. We demonstrate it for the case of Dirichlet conditions \eqref{4.12}.

Here, there are four ghost points which are computed using extrapolation of two points and the boundaries as in \eqref{5.10}, with additional terms:
\begin{eqnarray*} \label{7.12}
u_{-1/4} &=& -u_{1/4} + 2g_0 + \left( \frac{h}{4}\right)^2 u_{xx}(0,t) + \frac{1}{12}\left( \frac{h}{4}\right)^4 \frac{\partial^4 u}{\partial x^4}(0,t)+O(h^6) \nonumber \\
u_{-3/4} &=& -u_{3/4} + 2g_0 - 9\left( \frac{h}{4}\right)^2 u_{xx}(0,t) -\frac{81}{12} \left( \frac{h}{4}\right)^4 \frac{\partial^4 u}{\partial x^4}(0,t)+O(h^6) \nonumber \\
u_{N+1/4} &=&  -u_{N-1/4} + 2g_\pi + \left( \frac{h}{4}\right)^2 u_{xx}(\pi,t) + \frac{1}{12}\left( \frac{h}{4}\right)^4 \frac{\partial ^4 u}{\partial x^4}(\pi,t)+O(h^6) \nonumber \\
u_{N+3/4} &=& -u_{N-3/4} + 2g_\pi - 9\left( \frac{h}{4}\right) ^2 u_{xx}(\pi,t)-\frac{81}{12} \left( \frac{h}{4}\right)^4 \frac{\partial^4 u}{\partial x^4}(\pi,t) + O(h^6) 
\end{eqnarray*}
where $u_{xx}(0,t),u_{xx}(\pi,t)$ are expressed using the PDE as before \eqref{5.14} and
\begin{eqnarray*} \label{7.14} 
\frac{\partial^4 u}{\partial x^4}(0,t) &=& g_{0}''(t) - F_{t}(0,t) - F_{xx}(0,t) \\
\frac{\partial^4 u}{\partial x^4}(\pi,t) &=& g_{\pi}''(t) - F_{t}(\pi ,t)-F_{xx}(\pi,t) \ .
\end{eqnarray*}
%
%
Similarly, in the case of  Neumann conditions, the ghost points are:
\begin{eqnarray*}  \label{7.26}
u_{-1/4} &=& u_{1/4} - \left( \frac{h}{2}\right)u_{x}(0,t) - \frac{1}{3}\left( \frac{h}{4}\right) ^3 u_{xxx}(0,t) \nonumber \\
&& \hspace{3em}  \, - \, \frac{1}{60}\left( \frac{h}{4}\right) ^5 \frac{\partial^5 u}{\partial x^5}(0,t) + O(h^7) 
 \nonumber \\
u_{-3/4} &=& u_{3/4} - 3\left( \frac{h}{2}\right)u_{x}(0,t) - 9\left( \frac{h}{4}\right) ^3 u_{xxx}(0,t)  \nonumber \\
&& \hspace{3em} \, - \, \frac{243}{60}\left( \frac{h}{4}\right) ^5 \frac{\partial^5 u}{\partial x^5}(0,t) + O(h^7) \\
u_{N+1/4} &=& u_{N-1/4} + \left( \frac{h}{2}\right)u_{x}(\pi,t) + \frac{1}{3}\left( \frac{h}{4}\right)^3 u_{xxx}(\pi,t)  \nonumber \\
&& \hspace{3em} \, + \,  \frac{1}{60}\left( \frac{h}{4}\right) ^5 \frac{\partial^5 u}{\partial x^5}(\pi,t) + O(h^7) \\
u_{N+3/4} &=& u_{N-3/4} + 3\left( \frac{h}{2}\right)u_{x}(\pi,t) + 9\left( \frac{h}{4}\right) ^3 u_{xxx}(\pi,t)  \nonumber \\
&& \hspace{3em}  \, + \,
\frac{243}{60}\left( \frac{h}{5}\right) ^5 \frac{\partial^5 u}{\partial x^5}(\pi,t) + O(h^7)
\end{eqnarray*}
%
%
where $u_{xxx}(0,t),u_{xxx}(\pi,t)$ are expressed using the PDE as in \eqref{5.42} and
\begin{eqnarray*} \label{7.34} 
\frac{\partial^5 u}{\partial x^5}(0,t) &=& g_{0}''(t) - F_{tx}(0,t) - F_{xxx}(0,t) \\
\frac{\partial^5 u}{\partial x^5}(\pi,t) &=& g_{\pi}''(t) - F_{tx}(\pi,t) - F_{xxx}(\pi,t) \ .
\end{eqnarray*}
%
%

Consider the fourth order approximation of two-point block, as in \eqref{3.56}, for $j=1,...,N-2$:  
\begin{eqnarray*} \label{7.2}
\frac{d^2}{dx^2}u_{j+1/4} & \approx & \frac{1}{12(h/2)^2} [(-u_{j-3/4} + 16u_{j-1/4} - 30u_{j+1/4} + 16u_{j+3/4} - u_{j+5/4}) \\
&+& c(-u_{j-3/4} + 5u_{j-1/4} - 10u_{j+1/4} + 10u_{j+3/4} - 5u_{j+5/4} + u_{j+7/4})] \\
\frac{d^2}{dx^2}u_{j+3/4} & \approx & \frac{1}{12(h/2)^2} [(-u_{j-1/4} + 16u_{j+1/4} - 30u_{j+3/4} + 16u_{j+5/4} - u_{j+7/4}) \\
&+& c(u_{j-3/4} - 5u_{j-1/4} + 10u_{j+1/4} - 10u_{j+3/4} + 5u_{j+5/4} - u_{j+7/4})] 
\end{eqnarray*}
whereas near the boundaries we have:
\begin{eqnarray*}  \label{7.18}
\frac{d^2}{dx^2}u_{1/4} & \approx & \frac{1}{12(h/2)^2}   \Bigg[  (30+8c)g_0 + (7-4c)\left( \frac{h}{4}\right)^2 u_{xx}(0,t) \nonumber \\
&&  \hspace{2em} \, - \,  \frac{(65+76c)}{12} \left( \frac{h}{4}\right)^4 \frac{\partial^4 u}{\partial x^4}(0,t) 
  + (-46u_{1/4} + 17u_{3/4} - u_{5/4})  \nonumber \\
&&   \hspace{4em} \, + \, c(-15u_{1/4} + 11u_{3/4} - 5u_{5/4} + u_{7/4})  \Bigg]   \nonumber \\
\frac{d^2}{dx^2}u_{3/4} & \approx & \frac{1}{12(h/2)^2}   \Bigg[  (-2-8c)g_0 + (-1+4c)\left( \frac{h}{4}\right)^2 u_{xx}(0,t) +  \nonumber \\
&&  \hspace{2em}  \frac{(-1+76c)}{12} \left( \frac{h}{4}\right)^4 \frac{\partial ^4 u}{\partial x^4}(0,t) \nonumber \\
&&  \hspace{2em} \, + \, (17u_{1/4} - 30u_{3/4} + 16u_{5/4} - u_{7/4}) \nonumber \\
&&  \hspace{2em}  \, + \, c(15u_{1/4} - 11u_{3/4} + 5u_{5/4} - u_{7/4})  \Bigg]  \\
\frac{d^2}{dx^2}u_{N-3/4} & \approx & \frac{1}{12(h/2)^2}   \Bigg[  (-2-8c)g_{\pi} + (-1+4c)\left( \frac{h}{4}\right)^2 u_{xx}(\pi,t)  \nonumber \\
&&  \hspace{2em}  \, + \, \frac{(-1+76c)}{12} \left( \frac{h}{4}\right)^4 \frac{\partial^4 u}{\partial x^4}(\pi,t) \nonumber \\
&&  \hspace{2em} \, + \, (17u_{N-1/4} - 30u_{N-3/4} + 16u_{N-5/4} - u_{N-7/4})  \nonumber \\
&&  \hspace{2em} \, + \,  c(15u_{N-1/4} - 11u_{N-3/4} + 5u_{N-5/4} - u_{7/4})  \Bigg]  \\ 
\frac{d^2}{dx^2}u_{N-1/4} & \approx & \frac{1}{12(h/2)^2}   \Bigg[  (30+8c)g_{\pi} + (7-4c)\left( \frac{h}{4}\right)^2 u_{xx}(\pi,t) \nonumber \\
&&  \hspace{2em} \, - \, \frac{(65+76c)}{12} \left( \frac{h}{4}\right)^4 \frac{\partial ^4 u}{\partial x^4}(\pi,t) \nonumber \\
&&  \hspace{2em} \, + \, (-46u_{N-1/4} + 17u_{N-3/4} - u_{N-5/4}) 
 \nonumber \\
&&  \hspace{2em} \, + \, c(-15u_{N-1/4} + 11u_{N-3/4} - 5u_{N-5/4} + u_{N-7/4})  \Bigg] 
\end{eqnarray*}
%
This scheme becomes of fifth order for $c=\frac{4}{13}$.

\subsection{Numerical Example}
We run the scheme for the problem \eqref{4.10} with Dirichlet boundary conditions \eqref{4.12} with the same conditions, for $N=32,64,128,256,512$ grid points and fourth explicit Runge-Kutta time propagator.
In Figure \ref{fig:grid35} below, we compare the scheme for different values of $c$ and it can be clearly seen that for the value $c=\frac{4}{13}$ the scheme indeed becomes of fifth order at time $t=\pi$.
\begin{figure}[!htb]
  \begin{center}
    \includegraphics[clip,scale=0.5]{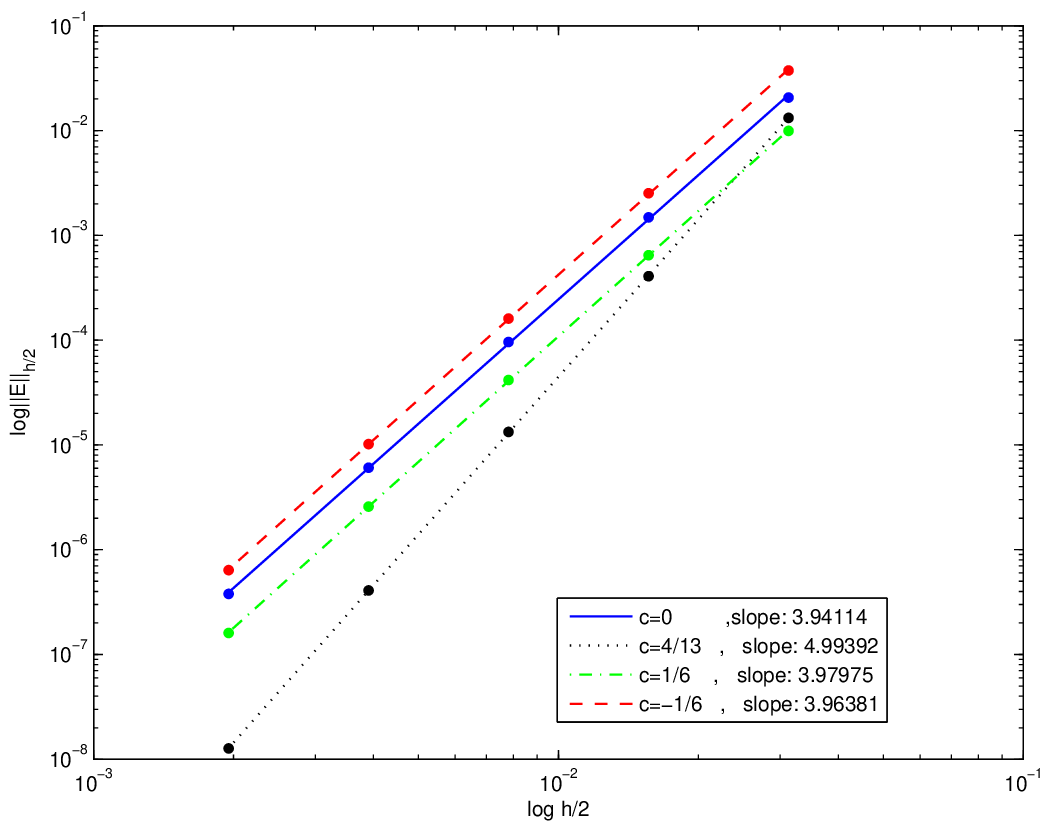}
            \caption{Convergence plot of fifth Order scheme for non-homogeneous Dirichlet \eqref{4.12}, $\log_{10} \Vert \textbf{E} \Vert \; vs. \; \log_{10} \left(  \frac{h}{2} \right)  $ for c=0, 4/13, 1/6, -1/6.}\label{fig:grid35}
  \end{center}
\end{figure}

\section{Conclusions}

This paper presented a novel methodology for designing BFD schemes whose global error is of a higher order than the truncation error. We considered the heat equation with periodic, Dirichlet, or Neumann boundary conditions. In each case, BFD schemes of third and fifth-order were derived by correctly defining the boundary stencils, such that the truncation errors lie in a different subspace than the solution, and construct the semi-discrete operators that guarantee error inhibition in time.

For the third-order schemes, we presented a thorough error and stability analysis. We presented numerical examples that demonstrate the theoretical results. 

We believe that our approach can be extended to other PDEs, such as advection equations,  and that even higher-order schemes can be obtained by using post-processing. We also expect that these ideas can be applied to other numerical methods, such as Finite Elements and Discontinuous Galerkin. These may be the topics for future research.

\bibliographystyle{elsart-num-sort}
\bibliography{myref2}





\end{document}